\documentclass[article,12pt,reqno]{amsart}
 \DeclareFontFamily{U}{wncy}{}
    \DeclareFontShape{U}{wncy}{m}{n}{<->wncyr10}{}
    \DeclareSymbolFont{mcy}{U}{wncy}{m}{n}
    \DeclareMathSymbol{\Sh}{\mathord}{mcy}{"58}
\usepackage{scalerel}
\usepackage[mathscr]{euscript}
\usepackage{hyperref}
\usepackage{amsmath}
\usepackage{amssymb}
\usepackage{amscd}
\usepackage{amsxtra}
\usepackage{mathrsfs}
\usepackage{algorithm}
\usepackage{algpseudocode}
\usepackage{tikz-cd}
\renewcommand{\AA}{\mathbb{A}}

\newcommand{\CC}{\mathbb{C}}

\newcommand{\FF}{\mathbb{F}}

\newcommand{\NN}{\mathbb{N}}
\newcommand{\PP}{\mathbb{P}}
\newcommand{\QQ}{\mathbb{Q}}
\newcommand{\RR}{\mathbb{R}}
\newcommand{\ZZ}{\mathbb{Z}}
\newcommand{\Oo}{\mathcal{O}}
\newcommand{\Cc}{\mathcal{C}}
\newcommand{\Dd}{\mathcal{D}}
\newcommand{\Ee}{\mathcal{E}}
\newcommand{\Hh}{\mathcal{H}}
\newcommand{\Gg}{\mathcal{G}}

\newcommand{\Ww}{\mathcal{W}}
\newcommand{\pr}{\mathfrak{p}}
\newcommand{\al}{\alpha}

\newcommand{\ep}{\varepsilon}

\newcommand{\si}{\sigma}

\newcommand{\inv}{^{-1}}
\newcommand{\set}[1]{\left\{ #1\right\}}
\newcommand{\ttm}[4]{\begin{pmatrix} #1 & #2 \\ #3 & #4 \end{pmatrix}}

\newcommand{\ip}[1]{\left\langle #1 \right\rangle }
\newcommand{\qr}[2]{\left( \frac{#1}{#2}\right)}

\newtheorem{thm}{Theorem}[section]

\newtheorem{lem}[thm]{Lemma}
\newtheorem{prp}[thm]{Proposition}
\newtheorem{cor}[thm]{Corollary}

\begin{document}

\title{Supersingular Isogeny Graphs from Algebraic Modular Curves}
\author{Nadir Hajouji}

\begin{abstract}
\noindent
We describe and compare algorithms for computing supersingular isogeny graphs.
Along the way, we obtain a formula for the trace of the adjacency matrix of a general supersingular isogeny graph,
and we prove a conjecture recently posed by Nakaya.
\end{abstract}

\maketitle

\section{Introduction}

Let $p>3$ be a prime. 
An elliptic curve $E/\FF_p^{sep}$ is \emph{supersingular} if the multiplication-by-$p$ endomorphism on $E$ is purely inseparable,
and $j(E) \in \FF_{p^2}$.\footnote{Note that the first condition actually implies the second.}
We write $S_p \subset \FF_{p^2}$ to denote the set of $j$-invariants of supersingular elliptic curves
$E/\FF_p^{sep}$, and we refer to elements of $S_p$ as supersingular $j$-invariants.

For a prime $\ell \neq p$, and an elliptic curve $E/\FF_p^{sep}$,
we write $S_\ell(E)$ to denote the set of subgroups $H \subset E(\FF_p^{sep})$ of order $\ell$.
For each $H \in S_\ell(E)$,
there is an associated isogeny $E \to E'$, with kernel equal to $H$, to some elliptic curve $E'$ determined by $H$.
Note that $E'$ is supersingular if and only if $E$ is supersingular.

Given a pair of distinct primes $p,\ell$,
we define the supersingular
isogeny graph $\Gamma_{p,\ell}$ as follows:

\begin{enumerate}
    \item The set of vertices of $\Gamma_{p,\ell}$ is the set $S_p$ of supersingular $j$-invariants.
    \item For each $j_0 \in S_p$, we fix a model $E_0/\FF_p^{sep}$ with $j(E_0) = j_0$.
    The graph has an edge $j_0 \to j(E_0/H)$ for each of the $\ell+1$ subgroups $H\in S_\ell(E_0)$.
\end{enumerate}

Our goal in this paper is to describe algorithms for computing $\Gamma_{p,\ell}$.

\section{General Background}

\subsection{General results}

We start by proving some general results.

\subsubsection{Frobenius}

In order to determine the field of definition of the objects needed to compute $\Gamma_{p,\ell}$,
we need to understand how Frobenius acts on supersingular curves.

{\lem{
\label{lem:frobss}
Let $p$ be a prime, let $q = p^{2n}$ let $E/\FF_q$ be a supersingular elliptic curve
and let $\phi_{2n} \in End(E)$ be the $q$-power Frobenius.

\begin{enumerate}
    \item $\phi_{2n} = \al \circ [p^n]$ for some $\al \in Aut(E)$.
    \item Suppose $\phi_{2n} = [\pm p^n]$ (i.e. $\al \in \set{[+1],[-1]}$),
    and let $N = p^n \mp 1$.
    \begin{enumerate}
        \item  $E(\FF_q) \cong \ZZ/N\ZZ \times \ZZ/N\ZZ$.
        \item Let $E\inv /\FF_q$ be the quadratic twist of $E$.
        The $q$-power Frobenius on $E\inv$ acts as $[\mp p^n]$.
    \end{enumerate}

\end{enumerate}
}}

\begin{proof}
    Let $[p^n] \in End(E)$ be the multiplication-by-$p^n$ endorphism.
    Since $E$ is supersingular, this endomorphism is purely inseparable,
    so by Cor II.2.12 in \cite{sil1}, there exists an isogeny $\al : E \to E$ such that $[p^n] = \phi_{2n} \circ \al$.
    Note that $\al$ has degree 1 since $[p^n],\phi_{2n}$ have the same degree, so $\al \in Aut(E)$, proving (1).

    Now assume $\phi_{2n}$ acts on $E$ as $[\pm p^n]$ and let
    let $N = p^n \mp 1$.
    \begin{itemize}
        \item Case 1: $N = p^n-1, \phi_{2n} = [+p^n]$.
        Since $p^n = N+1$, it follows that $p^n \equiv 1 \pmod N$.
        \item Case 2: $N = p^n + 1, \phi_{2n} = [-p^n]$.
        Then $-N + 1 = -p^n$, so again $-p^n \equiv 1 \pmod N$.
    \end{itemize}
    In both cases, $\phi_{2n}$ acts as the identity on $E[N]$,
    so $E[N] \subset E(\FF_q)$.
    Observe that the trace of $\phi_{2n}$ is $\pm 2p^n$,
    so by the Weil-conjectures:
    \[\# E(\FF_q) = q+1 - (\pm 2p^n) = (p^n \mp 1)^2 = N^2 \]
    Thus, $E(\FF_q) = E[N]$, proving the first claim in (2).

    Finally, let $E\inv/\FF_q$ be the quadratic twist of $E$.\footnote{Elliptic curves over finite fields of odd characteristc have exactly one quadratic twist, up to $\FF_q$-isomorphism.} 
    Then $\# E\inv (\FF_q) = 2(q+1) - \# E(\FF_q)$, so the trace of $\phi_{2n}$ is $\mp 2p^n$.
    Furthermore, $E\inv$ is also supersingular, so $\phi_{2n}'$ acts on $E\inv$
    as $\al' \circ [p^n]$ for some $\al \in Aut(E)$ by (1).

    To complete the proof, we just need to check a few cases:

    \begin{itemize}
        \item If $j(E) \not\in \set{0,1728}$, then $Aut(E) = Aut(E\inv) = \set{[\pm 1]}$,
        so $\phi_{2n}$ must act as $[+p^n]$ or $[-p^n]$.
        Only one of these possibilities has the correct trace.
        \item If $j(E) = 1728$ and $\al$ is the generator of the automorphism group,
        $\al \circ [p]$ has trace $0$.
        \item If $j(E) = 0$, $\al \in Aut(E)$ and $\al \neq [\pm 1]$, 
        then $\al \circ [p]$ has trace $\pm p$.
    \end{itemize}

    This completes the proof of (2).    
    
\end{proof}

Note that the assumption $\phi_{2n} = [\pm p^n]$ in (2) is automatically satisfied
when either of the following conditions hold:
\begin{itemize}
    \item If $E$ is defined over $\FF_p$:, then $\phi_2$ acts as $[-p]$ and $\phi_{2n}$ acts as $[(-p)^n]$.
    \item If $j(E) \not\in \set{0,1728}$, then $Aut(E) = \set{[\pm 1]}$,
    so $\phi_{2n} = [\pm p^n]$ by (1).
\end{itemize}

Thus, in all cases:
\begin{itemize}
    \item We can obtain a model for every $j\in S_p$ where $\phi_{2n}$ acts as $[\pm p^n]$.
    \item We can obtain a model where $\phi_{2n}$ acts as $[\mp p^n]$ by quadratically twisting.
\end{itemize}

\subsubsection{Induced Endomorphisms}

Let $\psi : E \to E'$ be a separable isogeny,
and define:
\[ End(E)_\psi = \set{ \xi \in End(E) : \ker(\psi) \subset \ker(\psi \circ \xi)}\]
For each $\xi \in End(E)_\psi$,
there exists a unique endomorphism $\xi_\psi \in End(E')$ that makes the obvious diagram commute:

\begin{center}
\begin{tikzcd}
E \arrow[r, "\psi"] \arrow[d, "\xi"] \arrow[d] & E' \arrow[d, "\exists ! \xi_\psi", dashed] \\
E \arrow[r, "\psi"]                            & E'                                        
\end{tikzcd}
\end{center}

{\lem{ \label{lem:leakpropscomp}
Let $\psi : E \to E'$ be a separable isogeny.
\begin{enumerate}
    \item $End(E)_\psi$ contains the multiplication-by-$m$ maps for all $m\in \ZZ$.
    \item The map $\xi\mapsto \xi_\psi$ preserves degree.
    \item For all $\xi_1, \xi_2 \in End(E)_\psi$, we have $\xi_1 \circ \xi_2 \in End(E)$ and $(\xi_1 \circ \xi_2)_\psi = (\xi_1)_\psi \circ (\xi_2)_\psi$.

\end{enumerate}

}}

\begin{proof}
    The first point follows from the definition of $End(E)_\psi$.

    Since $\psi \circ \xi = \xi_\psi \circ \hat{\psi}$ and $\psi, \hat{\psi}$ have the same degree, it follows that $\xi,\xi_\psi$ have the same degree, which proves the second point.

    If $\xi_1, \xi_2 \in End(E)_\psi$, we have $\ker(\psi) \subset \ker(\xi_2) \subset \ker(\xi_1 \circ \xi_2)$ so $\xi_1 \circ \xi_2$.
   Finally, we note that $(\xi_1 \circ \xi_2)_\psi = (\xi_1)_\psi \circ (\xi_2)_\psi$ follows from the fact that the right hand side satisfies the universal property defining the left hand side.
\end{proof}

Next, define:

{\lem{ \label{lem:leakpropsinj}

Let $\psi : E\to E'$ be a separable isogeny and $\hat{\psi} :E' \to E$ the dual isogeny.
Define a function $F_\psi : End(E)_\psi \to End(E')$ by $F_\psi(\xi) = \hat{\xi_\psi}$,
and define $F_{\hat{\psi}}$ analogously.
For all $\xi \in End(E)_\psi$, the following properties hold:
\begin{enumerate}
    \item $F_\psi(\xi) \in End(E')_{\hat{\psi}}$.
    \item $F_{\hat{\psi}}(F_\psi(\xi)) = \xi$, i.e. the maps $F_\psi,F_{\hat{\psi}}$ are inverses of each other. In particular, $F_{\hat{\psi}}$ is injective.

\end{enumerate}

}}
\begin{proof}
    By definition, we have $\psi \circ \xi = \xi_\psi \circ \psi$.
    Taking duals of both sides shows $\hat{\xi} \circ \hat{\psi} = \hat{\psi} \circ \hat{\xi_\psi} = \hat{\psi} \circ F_\psi(\xi)$.

    Now, $\ker(\hat{\psi}) \subset \ker(\hat{\xi} \circ \hat{\psi}) = \ker(\hat{\psi} \circ F_\psi(\xi))$, so $F_\psi(\xi) \in  End(E')_{\hat{\psi}}$, which proves (1).
    
    Furthermore, this shows that $(F_\psi(\xi))_{\hat{\psi}} = \hat{\xi}$ so $F_{\hat{\psi}} ((F_\psi(\xi))_{\hat{\psi}}) = \hat{\hat{\xi}} = \xi$. This holds for all $\xi \in End(E)_\psi$ and by symmetry, the same is true for all $\xi' \in End(E')_{\hat{\psi}}$, proving (2).
\end{proof}

We will return to this topic in the next section.

\subsection{Computing with models}
\label{ssec:modelstheory}

The standard way of computing $\Gamma_{p,\ell}$ is by doing the following:
\begin{itemize}
    \item Find a model representing every supersingular $j$-invariant.
    \item Compute $S_\ell(E)$ for each model.
    \item Use V{\'e}lu's formula (\cite{velu}) formulae to determine $j(E/H)$ for every $H \in S_\ell(E)$.
\end{itemize}

In order to use this approach, we need to find models of each supersingular curve with full $\ell$-torsion.
This may require working over a large field extension $\FF_q$,
but the next proposition shows that we already know how to obtain
the smallest such extension:

{\prp{\label{prp:models}

Let $p, \ell$ be primes, with $p\neq \ell$,
let $n$ be the order of $p^2$ in $(\ZZ/\ell \ZZ)^\times$
and let $q = p^{2n}$.

\begin{enumerate}
    \item For each $j \in S_p$, there is a model over $E/\FF_q$ with $j(E) = j$ and $E[\ell]\subset E(\FF_q)$.
    \item $\FF_q$ is the smallest field extension of $\FF_{p^2}$ where (1) holds.
\end{enumerate}

}}

\begin{proof}
    Let $j\in S_p$ and choose a model $E_0/\FF_q$ where $\phi_2$ acts as $[\pm p^n]$.
    Then $E_0(\FF_q) \cong \ZZ/N\ZZ \times \ZZ/N\ZZ$, where $N = p^n \mp 1$.

    \begin{itemize}
        \item If $\ell | N$, then $E_0$ has the desired properties.
        \item Otherwise, $\ell \not| p^n \mp 1$ but $\ell | p^{2n}-1$,
        so $\ell | p^n \pm 1$.
        Thus, the quadratic twist of $E_0$ has the desired properties.
    \end{itemize}

    Now, suppose (1) is true over $\FF_{p^{2m}}$ for some integer $m$.
    Since $S_p\neq \emptyset$, this means there exists at least one supersingular curve
    $E/\FF_{p^{2m}}$ with full $\ell$-torsion.
    This is only possible if $\FF_{p^{2m}}$ contains a primitive $\ell$th root of unity,
    so $\ell | p^{2m}-1$, i.e. $p^{2m} \equiv 1 \pmod \ell$.
    This means $m$ is a multiple of the order of $p^2$, so $\FF_q$ is isomorphic to a subfield of $\FF_{p^{2m}}$.

\end{proof}

Note that if $\ell = 2,3$, then $p^2 \equiv 1 \pmod \ell$ for all $p\neq \ell$,
so we can always compute $\Gamma_{p,\ell}$ over $\FF_{p^2}$ using this algorithm.

Now, to use this approach, we need a way of computing $S_\ell(E)$.

Prop \ref{prp:models} shows that, starting from any model $E_0/\FF_{p^2}$:
\[ y^2 = x^3 + fx + g\]
we can base change and/or quadratically twist to obtain a model given by an equation of the same form,
and $E(\FF_q) \cong (\ZZ/N) \times (\ZZ/N)$
for $N$ divisible by $\ell$.
Given such a model, we can compute $S_\ell(E)$ with the following algorithm:

\begin{algorithm}[ht!]
\label{alg:basistorsionsupsing}
	\caption{Computing $S_\ell(E)$ from model}
	\begin{algorithmic}
		\State{Set $ \texttt{Subgroups}= \set{}$.}
		\While{\texttt{True}}
			\State{Set $x$ to be a random element of $\FF_q$.}
			\State{Evaluate $y_2:=x^3 + fx + g$.}
			\If{$(\exists y\in \FF_q)(y^2=^3+fx=g)$:}
				\State Set $P = (x,y)$.
				\State Compute $Q = \frac{N}{\ell} \cdot P$.
				\If{$Q \neq O$}:
					\If{$\mathrm{length}(\texttt{Subgroups}) = 0$}
						\State Add $H = \set{n \cdot Q : 0\leq n < \ell}$ to the set \texttt{Subgroups}.
					\ElsIf{$\texttt{Subgroups} =\set{ H} $ and $Q \not\in H$}
						\For{$Q_i \in H$:}
							\State{Add $\set{n \cdot (Q+Q_i) : 0\leq n < \ell}$ to the set \texttt{Subgroups}.}
						\EndFor
						\State \Return \texttt{Subgroups}
					\EndIf
				\EndIf
			\EndIf
		\EndWhile
	\end{algorithmic}
\end{algorithm}

\begin{itemize}
    \item The probability that a random $x$ gives us a point on $E$ is $\frac{1}{2}$.
    \item If we find a point $P$ on $E$, the probability that $Q = \frac{N}{\ell} P$ is a nonzero point is $1-\frac{1}{\ell^2}$.
    \item If we already have a subgroup $H$ and we find a new point $P$,
    the probability that $P$ is not in that subgroup is $1-\frac{1}{\ell}$.
\end{itemize}

The expected number of trials before we obtain our first basis vector is $\frac{2\ell^2}{\ell^2-1}$ and $\frac{2\ell}{\ell - 1}$ for the second basis vector.
For $\ell = 2$, the expected number of trials for both steps is under 7,
and for $\ell > 2$ the expected number of trials is less than 6. 

The only challenge when using this approach is that it requires a way of efficiently computing square roots in $\FF_{p^{2n}}$.
If $n = 1$, we can work over $\FF_{p^2}$.
We can show that computing square roots in $\FF_{p^2}$ is no more difficult
than computing square roots in $\FF_p$.

{\lem{Assume we have the ability to compute square roots in $\FF_p$.
Fix a nonsquare $-d \in \FF_p$ and let $\rho = x + y \sqrt{-d} \in \FF_{p^2}$.

If we know $\rho^2 = a+b \sqrt{-d}$,
then we can solve for $\rho$.}}

\begin{proof}
\begin{enumerate}
To obtain $\rho$, we need to solve the system of equations:
\[\left\{ \begin{matrix}
x^2 - dy^2 &= a \\
2xy &= b
\end{matrix} \right. \]
	\item First, note that $y = 0$ if and only if $\rho^2 = x^2 = a$ is a square in $\FF_p$.
	Since we can compute square roots of elements in $\FF_p$, we may assume $y \neq 0$.
	
	\item The second equation allows us to solve for $x = \frac{b}{2y}$.
	We plug this in to the first equation to obtain the following quartic in $y$:
	\[ 4dy^4 + 4ay^2 - b^2 = 0\]
	\item Consider the following related quadratic polynomial:
	\[ 4dY^2 + 4a Y - b^2 = 0\]
	\begin{itemize}
		\item The discriminant of the quadratic is $(4a)^2-4(4d)(-b^2) = 16(a^2+db^2) = 16 N(a+b\sqrt{-d})$
		Note that this must be a square - in fact the square root is $4N(\rho) = 4(a^2+db^2)$.
		Since we have the ability to compute square roots in $\FF_p$, we can find a square root of $4N(\rho)$
		and thus solve the quadratic over $\FF_p$.
		\item 
		The product of the roots of the quadratic is equal to $\frac{-b^2}{4d}$.
		Now, this product is equal to $(-d)$ multiplied by a square, so the product of the roots is a nonsquare.
		This means that \emph{exactly} one of the roots of the quadratic is itself a square in $\FF_p$.
	\end{itemize}
	\item Now, if $Y_0$ is a root of the quadratic, then $y = \sqrt{Y_0}$ is a root of the quartic.
	Thus, we solve the quadratic, find the root which has a square root in $\FF_p$ and compute that square root to obtain $y$.
	We then compute $x = \frac{b}{2y}$ to obtain $\rho$.
\end{enumerate}

\end{proof}

This works well for certain pairs $p, \ell$,
but using these methods for $p = 79, \ell = 71$ would require the ability
to compute square roots in the field with $79^{70} \sim 2^{441}$ elements.

\subsection{Modular curves}

Next, we explain how to compute $\Gamma_{p,\ell}$ using modular curves.

Let $X_0(\ell)$ be the fine moduli space that represents isomorphism classes of pairs $(E,H)$,
where $E$ is an elliptic curve $E$ and $H \subset E$ is a (cyclic) subgroup of order $\ell$.
We write $j_\ell : X_0(\ell) \to X(1)$ to denote the $j$-map, which takes $[(E,H)]$ to $j(E)$,
and we write $\mathrm{Fr}_\ell : X_0(\ell) \to X_0(\ell)$ to denote the Fricke involution, which swaps $[(E,H)]$ and $[(E/H, \ker(\widehat{E\to E/H}))]$.

Then $X_0(\ell)(\FF_p^{sep})$ contains all information about isogenies of degree $\ell$ between elliptic curves over $\FF_p^{sep}$.
Furthermore, the points that represent pairs $(E,H)$,
where $E$ is supersingular,
are all defined over $\FF_{p^2}$:
we can always find models $E/\FF_{p^2}$ for each isomorphism class of supersingular curve where $\phi_2$ acts as $[\pm p]$,
and if $H \subset E$ is a subgroup of order $\ell$, then $\phi_2(H) = [\pm p] H = H$.
This shows that the pair $(E,H)$ is fixed by $\phi_2$,
so it is represented by a point in $X_0(\ell)(\FF_{p^2})$.

Now, given $P \in X_0(\ell)$, we can compute $j_\ell(P), j_\ell \circ \mathrm{Fr}_\ell(P)$
to determine the domain and codomain of the isogeny $E\to E/H$ associated to the pair $(E,H)$.
Thus, if we have a model of $X_0(\ell)$, together with formulae for $j_\ell, \mathrm{Fr}_\ell$,
we should in principle be able to compute $\Gamma_{p,\ell}$.

There are two issues we will have to address in order to use this method:
\begin{enumerate}
    \item First, we need a way of obtaining the models.
    For small values of $\ell$, these can be found in the literature.
    They can also be computed by hand, although this gets difficult as $\ell$ gets large.

    Fortunately, there is also a second modular curve we can use - the quotient of $X_0(\ell)$ by the Fricke involution. Models of the latter are easier to obtain.

    \item Once we have the model and the necessary formulae, we have to deal with the fact that \emph{some} points on $X_0(\ell)$ (or the quotient modular curve) may represent multiple edges on the isogeny graph. In order to obtain the correct graph, we need a way to determining the exact number of subgroups represented by a point on $X_0(\ell)$.
\end{enumerate}

The second problem is easier to work around.

Define $\Gamma_{p,\ell}^*$ to be the graph obtained from $\Gamma_{p,\ell}$
by removing all edges of the form $j_0 \to j_0$.
We will use modular curves to compute $\Gamma_{p,\ell}^*$,
and then obtain $\Gamma_{p,\ell}$ from $\Gamma_{p,\ell}^*$
by adding edges $j_0 \to j_0$ until each vertex has degree $\ell+1$.

\subsubsection{$X_0(\ell)$}

We start by explaining how to compute $\Gamma_{p,\ell}^*$
from the modular curve $X_0(\ell)$.
We require the following:

\begin{itemize}
    \item The set $S_p$ of supersingular $j$-invariants.
    \item A model of $X_0(\ell)$, together with formulae for $j_\ell,  \mathrm{Fr}_\ell$.
\end{itemize}

To compute the graph, we will do the following:
\begin{enumerate}
    \item Compute $j_\ell\inv(S_p)$.
    \item For each $P \in j_\ell\inv (S_p)$, we compute $j_\ell(P), j_\ell \circ \mathrm{Fr}_\ell(P)$.
    If the $j$-invariants are not equal, then we add edges $j_\ell(P) \to j_\ell \mathrm{Fr}_\ell(P) $ to the graph.
    The number of edges we add will depend on whether $j_\ell(P) \in \set{0,1728}$.

    \item Once we've checked every element in $j_\ell\inv(S_p)$, we add the missing edges to obtain $\Gamma_{p,\ell}$.
\end{enumerate}

{\lem{

Let $E$ be an elliptic curve and assume $j(E) \not\in\set{0,1728}$.
The function $S_\ell(E) \to X_0(\ell)$ that takes $(E,H)$ to the isomorphism class of $[(E,H)]$
is injective.

}}

\begin{proof}
    Let $H_1, H_2 \in S_\ell(E)$ and suppose $[(E,H_1)] = [(E,H_2)]$.
    Then there exists an automorphism $\al : E \to E$ satisfying $\al(H_1)= H_2$.
    Since $j(E) \not\in \set{0,1728}$, $Aut(E) = \set{[\pm 1]}$, so that means $\pm H_1 = H_2$.
    But $\pm H_1 = H_1$, so in fact we have $H_1 = H_2$.
    Thus, the map is injective.
\end{proof}

On the other hand, if $j(E) \in \set{0,1728}$,
the map will not be injective because $Aut(E)$ acts nontrivially on $S_\ell(E)$.
This means that points $P \in j_\ell\inv(\set{),1728} \cap S_p)$ may represent multiple edges
on the isogeny graph.

{\lem{\label{lem:fixedsubgroups}
Let $E$ be an elliptic curve with $j(E) \in \set{0,1728}$,
let $H \in S_\ell(E)$ and assume $H$ is fixed by $Aut(E)$.
Then $E/H \cong E$.

}}

\begin{proof}
    Let $\psi : E \to E/H$ be the quotient map.
    For any $\al \in Aut(E)$, we have:
  \[\ker(\psi \circ \al) = \al\inv(\ker(\psi)) = \al\inv(H) = H = \ker(\psi) \]
    so $Aut(E) \subset End(E)_\psi$.
    The function $F_\psi: End(E)_\psi \to End(E/H)$ is injective and preserves composition,
    so $End(E/H)$ contains a subgroup isomorphic to $Aut(E)$.
    Since $j(E) \in \set{0,1728}$, this is enough to deduce $E/H \cong E$.
\end{proof}

We can now show that $\Gamma_{p,\ell}^*$ can be computed from $X_0(\ell)$.
Furthermore, all of the computations happen over $\FF_{p^2}$:

{\prp{
Let $P \in X_0(\ell)$, and suppose $j_\ell(P) \in S_p$.

\begin{enumerate}
    \item $P \in X_0(\ell)(\FF_{p^2})$.
    \item Suppose $j_\ell(P) \neq j_\ell(\mathrm{Fr}_\ell(P))$.
    The number of edges $j_\ell(P) \to j_\ell(\mathrm{Fr}_\ell(P))$ is equal to $\frac{|Aut(E)|}{2}$, where $E$ is any elliptic curve satisfying $j(E) = j_\ell(P)$
\end{enumerate}

}}

\begin{proof}
    Fix a model $E/\FF_{p^2}$ with $j_\ell(E) = j_\ell(P)$,
    and where $\phi_2$ acts as $[\pm p]$.
    The point $P$ represents the isomorphism class of $(E,H)$ for some $H \in S_\ell(E)$.
    We compute:
    \[ \phi_2 \cdot P = [(\phi_2 \cdot E, \phi_2 \cdot H)] = [(E,[\pm p] H)] = [(E,H)] = P\]
    Thus, $P$ must be fixed by $\phi_2$ so $P \in X_0(\ell)(\FF_{p^2})$, proving (1).

    Now, $Aut(E)/\ip{[\pm 1]}$ acts on $S_\ell(E)$.
    Furthermore, if the orbit of $H$ has fewer than $|Aut(E)/\ip{[\pm 1]}|$ elements,
    then $j(E) = j(E/H) \in\set{0,1728}$ by Lemma \ref{lem:fixedsubgroups}.
    Thus, if $j(E) \neq j(E/H)$, then the orbit of $H$ contains precisely $Aut(E)/\ip{[\pm 1]}$ elements, proving (2).

\end{proof}

\subsubsection{$X_0(\ell)^+$}

Let $X_0(\ell)^+ = X_0(\ell)/\ip{\mathrm{Fr}_\ell}$ be the quotient of $X_0(\ell)$ by the Fricke involution.
For a point $Q = \set{P, \mathrm{Fr}_\ell(P)} \in X_0(\ell)^+$,
we define $a_\ell(Q) = j_\ell(P)+j_\ell \circ \mathrm{Fr}_\ell (P)$
and $b_\ell(Q) = j_\ell(P)\cdot j_\ell \circ \mathrm{Fr}_\ell (P)$.
We use these to construct a quadratic polynomial:
\[ R_\ell(x,Q) = x^2 - a_\ell(Q) x + b_\ell(Q)  = (x - j_\ell(P))(x-j_\ell\circ \mathrm{Fr}_\ell(P))\]
To recover the $j$-invariants of the curves represented by points above $Q$ in $X_0(\ell)$,
we simply evaluate $a_\ell(Q), b_\ell(Q)$ and solve the quadratic.

Now, there is a new ambiguity we have to deal with:
if $P$ is a fixed point of the Fricke involution,
then $R_\ell(x,\set{P,\mathrm{Fr}_\ell(P)}) = (x-j_\ell(P))^2$,
and the repeated root $j_\ell(P)$ gives rise to exactly one edge on the isogeny graph.
However, it is possible that $P \neq \mathrm{Fr}(\ell)(P)$,
but $R_\ell(x,Q) = (x-j_\ell(P))^2$ anyways.
In the second case, the point $Q$ represents two edges on the isogeny graph.

Fortunately, these ambiguities do not affect the computation of $\Gamma_{p,\ell}^*$:

{\prp{\label{prp:x0plusprops}

Let $Q \in X_0(\ell)^+$ and suppose $R_\ell(x,Q) = (x-j_1)(x-j_2)$ for $j_1, j_2 \in S_p$.
\begin{enumerate}
    \item $Q \in X_0(\ell)^+(\FF_{p^2})$.
    \item Assume $j_1 \neq j_2$. Then $Q$ gives rise to $m_1$ edges $j_1\to j_2$
    and $m_2$ edges $j_2\to j_1$, were $m_i = \frac{|Aut(E_i)|}{2}$.
\end{enumerate}

}}

\begin{proof}
    Let $P, \mathrm{Fr}_\ell(P)$ be the points over $Q$ in $X_0(\ell)$.
    Then $P, \mathrm{Fr}_\ell(P) \in X_0(\ell)(\FF_{p^2})$ and the quotient map is defined over $\ZZ$ so $Q \in X_0(\ell)(\FF_{p^2})$.

    If $j_1 \neq j_2$, then $j_\ell(P)\neq j_\ell \circ \mathrm{Fr}_\ell(P)$ so $P \neq \mathrm{Fr}_\ell(P)$.
    Thus, we have two points above $Q$ in $X_0(\ell)$,
    and the number of edges represented by each of those points is $\frac{|Aut(E)|}{2}$ by the proposition.

\end{proof}

\subsection{Supersingular Primes}

Let $p$ be a prime.
We say that $p$ is a \emph{supersingular prime} if one of the following conditions holds:
\begin{enumerate}
    \item $S_p \subset \FF_p$, i.e. every supersingular curve admits a model over $\FF_p$.
    \item $X_0(p)^+$ has genus 0.
    \item $p$ divides the order of the Monster group.
\end{enumerate}

The equivalence of these conditions is part of the theory of generalized Monstrous moonshine, see Section 9 of \cite{moonshine} .

For a supersingular prime $\ell$, the function field of $X_0(\ell)^+$ is generated
by a single element, and the functions $a_\ell, b_\ell$ can be described as polynomials in terms of that generator.
Furthermore, it is always possible to find a generator $y$ of the function field of $X_0(\ell)^+$ so that $a_\ell, b_\ell \in \ZZ[y]$ are monic polynomials of degree $\ell, \ell+1$, respectively.
This allows us to encode all of our data into a single polynomial $R_\ell(x,y) \in \ZZ[x,y]$.
These polynomials are known as Atkin modular polynomials,
and explicit formulae for $R_\ell(x,y)$ for all supersingular primes can be found in SageMath.

Let $A$ be a UFD and let $f(x), g(x) \in A[x]$.
    We write $\mathrm{Res}_x(f(x),g(x))$ to denote the resultant of the polynomials $f,g$ with respect to $x$. Note that the resultant is an element of $A$, and vanishes if and only if $f(x), g(x)$ have a common root.
    Furthermore, for $f(x) \in A[x]$, we write $\mathrm{rad}(f(x))$ to denote the radical of $f(x)$ - this is the unique squarefree polynomial with the same set of roots of $f(x)$.

For supersingular primes $\ell$, we can use resultants to extract several other polynomials of interest from the Atkin polynmomial $R_\ell(x,y)$.

\begin{enumerate}
    \item Let $F_\ell(x,j)$ be the classical modular polynomial of level $\ell$.\footnote{See Ex. 2.18,2.19 in \cite{sil2} or \cite{coxprimes}}.
    We can obtain $F_\ell$ from $R_\ell$ using the following formula:
    \[ F_\ell(x,j) = \frac{1}{(x-j)^{\ell+1}}\mathrm{Res}_y (R_\ell(x,y),R_\ell(j,y)) \]

    \item Define:
    \[ \Delta_\ell(x) = \mathrm{Res}_y\left( \frac{\partial}{\partial x} R_\ell(x,y), a_\ell(y)^2-4b_\ell(y)\right)\]
    The roots of this polynomial are precisely the points that represent endmorphisms of degree $\ell$ in $X_0(\ell)^+$.
    Note that $\Delta_\ell(x)$ can also be obtained (up to multiplication by a scalar) using the formula $\Delta_\ell(x) = F_\ell(x,x)$.
    
    We will describe this polynomial in detail in the next section.

    \item Let $s_p(x) = \prod_{j \in S_p} x- j$ and define the \emph{supersingular polynomial of level $\ell$}:
    \[ s_{p,\ell}^+(y) = \mathrm{rad} \left(\mathrm{Res}_x (R_\ell(x,y),s_p(x)) \right)\]

    The roots of this polynomial are precisely the points that represent isogenies between supersingular curves.
        Thus, the problem of finding all $Q \in X_0(\ell)^+$ that represent isogenies between supersingular curves is equivalent to computing the roots of $s_{p,\ell}^+$.

        The polynomials $s_{p,\ell}^+$ are also the subject of several conjectures 
        in \cite{x0plusconjecture}. 
        We end this section by using our results 
        to shed some light on questions raised in that paper.

\end{enumerate}

\subsection{Application: Nakaya's conjectures}
\label{ssec:nakaya}

We will focus on the questions and conjectures related to the splitting field of $s_{p,\ell}^+$.

To begin, we note that part (1) of Prop \ref{prp:x0plusprops} 
implies that $s_{p,\ell}^+$ always splits over $\FF_{p^2}$.
This explains why the answer to Question 1 in \cite{x0plusconjecture} is \emph{yes}.
On the other hand, we can use connectedness of $\Gamma_{p,\ell}$ to show that while $s_{p,\ell}^+$ splits over $\FF_{p^2}$,
it very rarely splits  over $\FF_p$:

{\prp{\label{prp:ssplnonsplit}
Let $\ell$ be a supersingular prime,
let $p\neq\ell$ be a prime and suppose $s_{p,\ell}^+$ splits over $\FF_p$.
Then $p$ is also a supersingular prime.

}}

\begin{proof}

Let $V_0 = S_p \cap \FF_p$ and $V_1 = S_p \backslash \FF_p$.
Then:
\begin{itemize}
    \item $V_0 \neq \emptyset$ for all primes $p$.
    \item $V_1 = \emptyset$ if and only if $p$ is supersingular.
\end{itemize}
For a prime $p$ which is \emph{not} supersingular,
the sets $V_0, V_1$ form a nontrivial partition of the set of vertices $S_p$
of $\Gamma_{p,\ell}$.
Since $\Gamma_{p,\ell}$ is connected,
it must contain an edge $j_0\to j_1$, where $j_0 \in V_0, j_1 \in V_1$.
This means there exists $y_0$ such that:
\[ R_\ell(x,y_0) = x^2 - a_\ell(y_0) x + b_\ell(y_0) =  (x - j_0)(x-j_1) \]
The point $y_0$ represents an isogeny between supersingular curves,
so it must be a root of $s_{p,\ell}^+(y)$.
Furthermore, since $a_\ell \in \ZZ[x]$, and $a_\ell(y_0) = j_0 + j_1 \not\in \FF_p$, it follows that $y_0 \not\in \FF_p$,
which proves that $s_{p,\ell}^+(x)$ does not split over $\FF_p$.

\end{proof}

Generalized moonshine allows us to associate a sporadic group $G_\ell$ to each of
the supersingular primes $\ell <12$,
and Conjecture 2 in \cite{x0plusconjecture} predicts that for $\ell \in \set{5,7}$,
$s_{p,\ell}^+$ splits over $\FF_p$ if and only if $p | |G_\ell|$.
Prop \ref{prp:ssplnonsplit} reduces the proof of this conjecture to a finite computation.

For a pair $(p,\ell)$, where $p,\ell$ are supersingular primes,
we write $\mathrm{Ob}(p,\ell)$ to denote the number of quadratic factors that appear in the factorization of $s_{p,\ell}^+$ over $\FF_p$.
We use this to quantify the \emph{obstruction} to $s_{p,\ell}$ splitting over $\FF_p$:
 $s_{p,\ell}^+$ splits over $\FF_p$ if and only if $\mathrm{Ob}(p,\ell) = 0$.
We computed the value of $\mathrm{Ob}(p,\ell)$ for all pairs $(p,\ell)$,
where $p,\ell$ are supersingular,
and recorded the values in Table \ref{table:ob}.

\begin{table}[ht]
\caption{$\mathrm{Ob}(p,\ell)$}
\label{table:ob}
\begin{tabular}{llllllllllllllll}
            & \textbf{2} & \textbf{3} & \textbf{5} & \textbf{7} & \textbf{11} & \textbf{13} & \textbf{17} & \textbf{19} & \textbf{23} & \textbf{29} & \textbf{31} & \textbf{41} & \textbf{47} & \textbf{59} & \textbf{71} \\
\textbf{2}  & 0          & 0          & 0          & 0          & 0           & 0           & 0           & 0           & 0           & 1           & 0           & 1           & 0           & 1           & 1           \\
\textbf{3}  & 0          & 0          & 0          & 0          & 0           & 0           & 0           & 1           & 0           & 0           & 2           & 1           & 1           & 2           & 2           \\
\textbf{5}  & 0          & 0          & 0          & 0          & 0           & 1           & 2           & 0           & 2           & 1           & 1           & 2           & 5           & 3           & 4           \\
\textbf{7}  & 0          & 0          & 0          & 0          & 1           & 1           & 0           & 2           & 2           & 3           & 3           & 2           & 3           & 5           & 9           \\
\textbf{11} & 0          & 0          & 0          & 1          & 0           & 1           & 3           & 2           & 5           & 4           & 4           & 6           & 10          & 12          & 15          \\
\textbf{13} & 0          & 0          & 1          & 1          & 1           & 0           & 2           & 4           & 2           & 5           & 8           & 11          & 9           & 12          & 16          \\
\textbf{17} & 0          & 0          & 2          & 0          & 3           & 2           & 0           & 5           & 5           & 10          & 7           & 15          & 12          & 18          & 22          \\
\textbf{19} & 0          & 1          & 0          & 2          & 2           & 4           & 5           & 0           & 7           & 5           & 10          & 11          & 16          & 17          & 20          \\
\textbf{23} & 0          & 0          & 2          & 2          & 5           & 2           & 5           & 7           & 0           & 12          & 12          & 16          & 21          & 27          & 33          \\
\textbf{29} & 1          & 0          & 1          & 3          & 4           & 5           & 10          & 5           & 12          & 0           & 12          & 24          & 26          & 28          & 38          \\
\textbf{31} & 0          & 2          & 1          & 3          & 4           & 8           & 7           & 10          & 12          & 12          & 0           & 17          & 28          & 34          & 39          \\
\textbf{41} & 1          & 1          & 2          & 2          & 6           & 11          & 15          & 11          & 16          & 24          & 17          & 0           & 36          & 46          & 51          \\
\textbf{47} & 0          & 1          & 5          & 3          & 10          & 9           & 12          & 16          & 21          & 26          & 28          & 36          & 0           & 56          & 67          \\
\textbf{59} & 1          & 2          & 3          & 5          & 12          & 12          & 18          & 17          & 27          & 28          & 34          & 46          & 56          & 0           & 83          \\
\textbf{71} & 1          & 2          & 4          & 9          & 15          & 16          & 22          & 20          & 33          & 38          & 39          & 51          & 67          & 83          & 0          
\end{tabular}
\end{table}

Note that Table \ref{table:ob}. has a property we did not expect it to have: it is \emph{symmetric}.
In other words, for all pairs of supersingular primes $p,\ell$, we have $\mathrm{Ob}(p,\ell) = \mathrm{Ob}(\ell,p)$.
This generalizes the "curious observation'' at the end of \cite{x0plusconjecture}.
Thus, one can think of the rows as representing either $p$ or $\ell$.

For concreteness, we will view the rows as representing the values of $\mathrm{Ob}(-,\ell)$.
Then we can read off the following:
\begin{itemize}
    \item The rows associated to $\ell = 5, 7$ confirm that Conjecture 2 is true.
    \item The sporadic group associated to $\ell = 11$ is the Mathieu group $M_{12}$.
    The primes that divide $M_{12}$ are $p = 2,3,5,11$,
    and the row $\ell = 11$ shows that these are precisely the primes where $\mathrm{Ob}(p,11) = 0$.
    
\end{itemize}

We will see that $\mathrm{Ob}(p,\ell)$ can sometimes be obtained without needing a formula for $R_\ell$ in the next section.

\section{Endomorphisms of degree $\ell$}
\subsection{Overview}

\subsubsection{Example}
\label{sssec:exp7l2}
Before introducing any new tools,
we investigate the special case where $\ell = 2$.
It is well-known that:
\[ \Delta_2(x) = (x-1728)(x-8000)(x+3375)^2\]
Observe that the roots of $\Delta_2$ are all congruent to $-1 \pmod 7$,
i.e.:
\[ \Delta_2(x) \equiv (x+1)^4 \pmod 7\]
This would suggest that the curve with $j = -1$ has 4 endomorphisms of degree 2,
but that is clearly impossible, 
since elliptic curves can't have more than 3 distinct subgroups of order 2.
We will show that this happens because the factor $(x+3375)^2$
behaves differently in characteristic 7.

\subsubsection{Notation}
\label{sssec:notation}

Let $-d <0$ be a fundamental discriminant and define:
\begin{itemize}
	\item $K = \QQ(\sqrt{-d})$.
	\item $\Oo_d$ is the ring of integers of $K$.
	\item $H_{-d}(x)$ is the Hilbert class polynomial of $K$,
            and $h(-d) = \deg H_{-d}(x)$.
	\item $L/K$ is the splitting field of $H_{-d}(x)$.
	\item $j_1, ..., j_h \in L$ are the roots of $H_{-d}(x)$, 
        and $E_1,\hdots, E_h/L$ are elliptic curves with $j(E_i) = j_i$.
	We write $[]_{E_i} : \Oo_d \to End(E_i)$ to denote the normalized                isomorphism described in II.1.1 of \cite{sil2}.
	\item $\Gg$ is the Galois group of $L/\QQ$.
        We write $\Hh$ to denote $Gal(L/K)$ and $\si$ to denote the generator of  $Gal(K/\QQ)$.
        Note that $\Hh$ is normal in $\Gg$,
        and $\Gg/\Hh \cong \ip{\si_d}$.
        
\end{itemize}

Let $m>1$ be an integer, let $X_0(m)$ be the modular curve of level $m$.
We define:
\begin{align*}
    X_0(m)_\Delta &= \set{P \in X_0(m): j_m(P) = j_m \circ \mathrm{Fr}_m(P)} \\
    X_0(m)[G]_\Delta = {P \in X_0(m)_\Delta : G(j_m(P)) = 0}
\end{align*}

Note that $X_0(m)[H_{-d}(x)]_\Delta = \emptyset$ unless $\Oo_d$ contains an element
of norm $m$.

Let $\ell$ be a prime.
We will write $\chi_\ell : \ZZ \to \NN$ to denote the trivial character:
\[ \chi_\ell(d) = \left\{
\begin{matrix}
    0 & \mathrm{if} & \ell | d \\
    1 & \mathrm{otherwise} & 
\end{matrix}
\right.\]

Now, assume $\ell$ splits in $\Oo_d$.\footnote{In other words, $\Oo_d$ contains an element of norm $\ell$.}]
We define an element $\xi_{\ell,-d}^+ \in \Oo_d$ as follows:
\begin{itemize}
    \item If $-d \not\in \set{-3,-4}$, and $\frac{a+b\sqrt{-d}}{2}$ has norm $\ell$, we set $\xi_{\ell,-d}^+ = \frac{|a|+|b|\sqrt{-d}}{2}$.
    \item If $-d = -4$, then $\ell = a^2+b^2$ for integers $a,b$.
    Let $a^+ = \min(|a|,|b|)$ and $b^+ = \max(|a|,|b|)$.
    We set $\xi_{\ell,-4}^+ = a^+\sqrt{-1} b^+$.
    \item If $-d = -3$, then $\ell = a^2 + 3b^2$ for some $a,b \in \ZZ$.
    We define $\xi_{\ell,-3}^+ = |a| + \sqrt{-3}|b|$.
\end{itemize}

Note that every element of norm $\ell$ can be obtained from $\xi_{\ell,-d}^+$
by multiplying by a unit in $\Oo_d^\times$, and/or conjugating by $\si_d$.

{\prp{\label{prp:frickechar0}
Suppose $\Oo_d$ contains an element of norm $\ell$,
where $\ell$ is prime,
and let $P \in X_0(\ell)[H_{-d}(x)]_\Delta$.

Then:
\begin{enumerate}
    \item If $P'\in  X_0(\ell)[H_{-d}(x)]_\Delta$ and $j_\ell(P') = j_\ell(P)$,
    then $P' = P$ or $P' = \si_d(P)$.
    
    \item $\si_d(P) = \mathrm{Fr}_\ell(P)$
    \item $\mathrm{Fr}_\ell(P) = P$ if and only if $\ell | d$.
\end{enumerate}

}}

\begin{proof}
    Let $E$ be an elliptic curve with CM, let $[\xi] \in End(E)$
    be an endomorphism of degree $\ell$,
    and let $\al \in Aut(E)$.
    Note that $\al \circ [\xi] = [\xi] \circ \al$,
    and $\ker(\al \circ [\xi]) = \ker([\xi])$ because $\al$ is injective.
    Thus, replacing $[\xi]$ by $\al \circ [\xi]$,
    we may assume that $\xi = \xi^+$ or $\xi = \si_d(\xi)^+$.

    In particular, this means that any single curve $E$ gives rise to at most 2 points in $X_0(\ell)_\Delta$, namely $P^+ = [(E,\ker([\xi^+]))]$ and $\si_d(P) = [(E,\ker(\si_d \cdot (\xi^+))]$.
    Furthermore, an easy computation shows that:
    \[[\xi^+] \circ [\si_d \xi^+] = [\xi^+ \circ \si_d \xi^+] = [N_{K/\QQ}(\xi^+] = [\deg (\xi)] \]
    so $[\si_d \cdot\xi^+] = \widehat{[\xi]}$.
    Finally, by II.1.1.1 in \cite{sil2},
    $[\si_d \cdot \xi] = \si_d \cdot [\xi]$,
    so $(E, \ker([\si_d \cdot\xi^+] )) = \si_d \cdot (E,\ker([\xi]))$.
    Thus:
    \[ \mathrm{Fr}_\ell((E,\ker([\xi]))) = [(E,\ker(\widehat{\xi}))] = \si_d \cdot [(E,\ker([\xi]))]) \]
    proving the first two claims.

    We now turn our attention to the third claim.
    If $\ell | d$, then $\ell$ ramifies in $\Oo_d$,
    so $\xi_{\ell,-d}^+,\si_d(\xi_{\ell,-d}^+)$ differ by a unit in $\Oo_d$.
    This means that the endomorphisms they represent have the same kernel,
    so they are represented by the same point on $X_0(\ell)$.
    Thus, $\mathrm{Fr}_\ell(P) = P$.

    Finally, we prove that $\mathrm{Fr}_\ell(P) = P \implies \ell | d$.
    Note that this can be verified by hand for $\ell = 2,3$,
    so we will assume that $\ell > 3$.

    The kernels of $\xi^+_{\ell,-d},\widehat{\xi^+_{\ell,-d}}$ are contained in $E[\ell]$, so we only need to study the restrictions of these endomorphisms to the torsion subgroup $E[\ell]$.
    We can view $E[\ell]$ as a 2-dimensional vector space over $\FF_\ell$,
    and the restrictions of the two endomorphisms as linear transformations.
    From this perspective, points in the kernel of $\xi_{p,\ell}^+$
    are the same as eigenvectors with eigenvalue 0.

    Let $a = Tr_{K/\QQ}(\xi^+_{\ell,-d}) \in \ZZ$.
    Then $\widehat{[\xi^+_{\ell,-d}]} = [\si_d \xi_{\ell,-d}^+] = [a] - [\xi_{\ell,-d}]$, so points in the kernel of $\widehat{[\xi^+_{\ell,-d}]}$
    are eigenvectors of  $\xi_{p,\ell}^+$ with eigenvalue $a$.
    In order to prove that the kernels are distinct,
    it suffices to show that $a \not\equiv 0 \pmod \ell$.

    Now, the minimal polynomial of $\xi^+_{\ell,-d}$ over $\QQ$ is $x^2- ax + \ell$.
    Since $\xi^+_{\ell,-d} \not \in \RR$,
    we have $a^2 < 4\ell$ and since $\ell \geq 5$,
    we deduce $a^2 < \ell^2$, so $|a|< \ell$.
    Thus $a \equiv 0 \pmod \ell$ if and only if $a = 0$,
    and $a = 0$ is equivalent to $\xi_{\ell,-d}^+ = \sqrt{-\ell}$.

\end{proof}

Let $P_0 \in X_0(\ell)[H_{-d}(x)]_\Delta$,
and let $\tau \in \Gg_d$.
If $P_0$ represents $(E_0, \ker([\xi]))$,
then $\tau \cdot P_0$ represents $(\tau \cdot E_0, \ker([\tau \cdot \xi]))$

For a general point $P_1 \in X_0(\ell)[H_{-d}(x)]$,
which represents $(E_1, H_1)$,
we can find $\tau \in \Hh$ such that $\tau E_0 = E_1$.
The lemma shows that either $\tau \cdot P_0 = P_1$,
or $\si_d \tau \cdot P_0 = P_1$.

We deduce:
\begin{itemize}
    \item $\Gg_d$ acts transitively on $ X_0(\ell)[H_{-d}(x)]_\Delta$.
    \item If $\ell \not| d$ and $\Oo_d$ contains an element of norm $\ell$,
    the action of $\Gg_d$ is simply transitive.
    \item If $\ell | d$ and $\Oo_d$ contains an element of norm $\ell$,
    the restriction of the action to $\Hh_d$ is simply transitive.
\end{itemize}

{\prp{\label{prp:cmsupsingredprops}

Let $-d$ be a discriminant and $p$ a prime satisfying $\qr{-d}{p}\neq 1$.

\begin{enumerate}
    \item For each root $j_0$ of $H_{-d}$,
    There exists $\tau \in \Hh$ satisfying:
    \begin{itemize}
	\item $j(\tau) = j_0$.
	\item $\tau = \frac{-b+\sqrt{-d}}{2a}$ for some $a\in \ZZ$, where $p\not| a$.
    \end{itemize}

    \item Let $j_1,j_2$ be roots of $H_{-d}$.
    There exists an isogeny $\psi : E_0\to E_1$ satisfying:
    \begin{itemize}
        \item $j(E_i) = j_i$ for $i =1,2$.
        \item The degree of $\psi$ is not divisible by $p$.
    \end{itemize}
\end{enumerate}

}}

\begin{proof}
    
    Let $q(x,y) = ax^2 + bxy + cy^2$ be a binary form with discriminant $b^2 - 4ac =-d$,
    and suppose $p | ac$.
    Note that this forces $-d = b^2 - 4ac \equiv b^2 \pmod p$.
    If $p\not| b$, then $-d$ has a nonzero square root in $\FF_p$ so $\qr{-d}{p} = 1$, which is a contradiction,
    so we deduce that $p|ac \implies p| b$.
    
    Now suppose $p | a$ and $p | c$.
    By the previous argument, $p | b$, so $p^2 | b^2 - 4ac = -d$.
    This contradicts the assumption that $-d$ is a fundamental discriminant,
    so we deduce that either $p\not| a $ or $p\not| c$.
    
    If $p \not| a$, we take $\tau$ to be the unique root of $ax^2 + bx + c = 0$ in the upper half plane.
    By the quadratic formula, $\tau = \frac{-b+\sqrt{-d}}{2a}$ as claimed in the lemma.
    Otherwise, we take $\tau = \frac{-b+\sqrt{-d}}{2c}$ to be the root of $cx^2 + bx + a = 0$ in the upper half plane.
    The two choices of $\tau$ lie in the same orbit under $z\mapsto \frac{-1}{z}$, so they lie in the same $SL_2(\ZZ)$-orbit of $\Hh$,
    and thus represent the same elliptic curve over $\CC$.
    \footnote{Note, however, that if $\tau$ represents an elliptic curve which is defined over $\RR$,
    then $\frac{-1}{\tau}$ represents the quadratic twist by $\sqrt{-1}$ of that elliptic curve.}
    This proves (1)

    We now turn our attention to (2).
    Define $\tau_0 = \frac{1+\sqrt{-d}}{2}$ if $-d \equiv 1 \pmod 4$ and $\tau_0 = \frac{\sqrt{-d}}{2}$ otherwise.
    Note that $\Oo_d = \ZZ \oplus \ZZ[\tau_0]$ is the ring of integers of $\QQ(\sqrt{-d})$.
    Furthermore, the elliptic curve $E_0 = \CC/\Oo_d$ has CM by $\Oo_d$.
    We set $j_0 = j(E_0)$ - note that $j_0$ is also a root of $H_{-d}(x)$.
    
    By (1),
    we can find $\tau_1, \tau_2 \in \Hh$,
    where $\tau_i = \frac{-b_i + \sqrt{-d}}{2a_i}$,
    such that $j(\tau_i) = j_i$ and $p\not| a_i$.
    We set $E_i = \CC /(\ZZ\oplus \tau_i \ZZ)$.

    Now, let $T_i = \tau_i + \Oo_d \in \CC/\Oo_d = E_0(\CC)$.
    Then $[a_i] \cdot T_i \in \Oo_d$, so $T_i$ is a point of order $a_i$.
    Furthermore, the quotient $E_0/\ip{T_i}$ is isomorphic to $E_i$.
    Thus, $E_0$ has isogenies $E_0 \to E_1$, $E_0\to E_2$ of degree $a_1,a_2$,
    respectively, and composing the dual $E_1\to E_0$ with $E_0 \to E_2$
    gives us an isogeny of degree $E_1 \to E_2$ of degree $a_1a_2 \neq 0 \pmod p$.
    
\end{proof}

Note that the condition $\qr{-d}{p}\neq 1$ is equivalent to the roots of $H_{-d}(x) \in \FF_p^{sep}$ being supersingular $j$-invariants by Ex.2.30 in \cite{sil2}.

The polynomial $\Delta_2$ is a product of Hilbert class polynomials,
so we have all of the tools we will need.
However, starting with $\ell = 3$, there are also factors that appear in $\Delta_\ell$ which are not Hilbert class polynomials, e.g.:
\begin{align*}
    \Delta_3(x) &= x(x-54000)(x-8000)^2(x+32768)^2\\
    &= H_{-3}(x)H_{-8}(x)^2H_{-11}(x)^2(x-54000)
\end{align*}

Now, $x-54000$ is not a Hilbert class polynomial,
but $j = 54000$ represents a curve in the same isogeny class as $j = 0$.
To analyze $\Delta_\ell$ for general $\ell$,
we will need names for these factors.

For a discriminant $d$ and a positive integer $m$,
we write $H_{-d,m}(x)$ to denote the (unique squarefree, monic) polynomial
whose roots are $j$-invariants of elliptic curves $E'$
for which there exists an isogeny $E\to E'$ of degree dividing $m$
from an elliptic curve $E$ satisfying $H_{-d}(j(E)) = 0$.
For example, $H_{-3,2}(x) = x(x-54000)$.

\begin{itemize}
    \item By definition, we have $H_{-d,1} = H_{-d}$; to keep the notation consistent, we will use $H_{-d,1}$ to refer to the Hilbert class polynomial from now on.
    This allows us to write $\Delta_3$ as:
    \[ \Delta_3(x) = H_{-3,2}(x)H_{-8,1}(x)^2H_{-11,1}(x)^2\]
    \item The roots of $H_{-d,m}$ always represent elliptic curves in the same isogeny class as $H_{-d}$, so we deduce that the roots of $H_{-d,m}$ are supersingular if and only if $\qr{-d}{p} \neq 1$.
\end{itemize}

Now, we can always find a factorization for $\Delta_\ell(x)$
in terms of $H_{-d,m}$ for some finite set of pairs $(-d,m)$.
In order to explain the presence of these factors $H_{-d,m}$ in $\Delta_\ell$,
and to determine which pairs $-d,m$ are needed, we need to extend some of the results from the previous section.

Let $\psi : E \to E"$ be a separable isogeny and define:
\[ \overline{End(E)_\psi} = \set{\xi \in End(E) : (\exists \xi' \in End(E'))(\psi \circ \xi = \xi' \circ \psi)} \]
Note that $End(E)_\psi \subset \overline{End(E)_\psi}$.

{\lem{
Let $\psi : E \to E'$ be a (nonzero) isogeny, let $\xi \in \overline{End(E)_\psi}$
let $\xi_1',\xi_2' \in End(E')$ and suppose $\psi \circ \xi = \xi_1' \circ \psi = \xi_2' \circ \psi$.
Then $\xi_1' = \xi_2'$.
}}

\begin{proof}
    If $(\xi_1 - \xi_2) \circ \psi = 0$,
    then we can compose with the dual of $\xi_1 - \xi_2$ to show that $[\deg (\xi_1 - \xi_2)] \circ \psi = 0$.
    Since $Hom(E, E')$ is a torsion-free $\ZZ$-module and $\psi$ is nonzero,
    it follows $\xi_1 = \xi_2$.
    
\end{proof}

Thus, the maps $\xi\mapsto \xi_\psi$ and $F_\psi: \xi \mapsto \widehat{\xi_\psi}$ are well-defined (and injective) on $\overline{End(E)}_\psi$,
not just $End(E)_\psi$.

{\lem{
For all $\xi \in \overline{End(E)}_\psi$ and $a \in \ZZ$,
we have $\xi +[a] \in \overline{End(E)}$.

}}
\begin{proof}
    This is obvious, since we can take $(\xi+[a]_E)_\psi = \xi_\psi + [a]_{E'}$.

\end{proof}

{\lem{
Let $E$ be an elliptic curve with CM by $\Oo_d$.
Fix $\xi_0 \in \Oo_d$ such that $\Oo_d = \ZZ[\xi_0]$,
and let $\xi = a + m \xi_0 \in \Oo_d$ be an element of norm $\ell$.

Let $\psi : E \to E'$ be an isogeny of degree dividing $m$.
Then $\xi \in \overline{End(E)_\psi}$.
}}

\begin{proof}
    First, note that $[m\xi_0] \in End(E)_\psi$,
    since $\ker(\psi) \subset E[m] \subset \ker(\psi \circ [m\xi_0])$.
    The result now follows from the previous lemma.

\end{proof}

If $\xi \in \Oo_d$ has norm $\ell$,
then $\xi$ satisfies a polynomial:
\[ x^2 + ax + \ell = 0\]
for some $a$.
Since $\Oo_d$ is an imaginary extension,
we must have $a^2 - 4\ell < 0$,
so there are only finitely many polynomials of this form for any fixed $\ell$.
Furthermore, $\xi\in\Oo_d$ implies that $a^2 - 4\ell = (-d)m^2$ for some integer $m$,
so $\xi = \frac{-a + m \sqrt{-d}}{2} = a' + m \xi_0$,
where $\xi_0 = \frac{\sqrt{-d}}{2}$ or $\xi_0 = \frac{1+\sqrt{-d}}{2}$.

Now, if $\psi: E\to E'$ is an isogeny of degree $m$ from an elliptic curve with CM by $\Oo_d$,
then $\ker(\psi) \subset E[m] \subset \ker([m] \circ \xi_0)$,
so $[m] \circ \xi_0 \in End(E)_\psi$,
and thus $\xi \in \overline{End(E)_\psi}$ by the lemma.
This shows that $E'$ also has an endomorphism of degree $\ell$,
so $\Delta_\ell(j(E')) = 0$, and more generally $H_{-d,m}(x) | \Delta_\ell(x)$.
Motivated by this, we write $\Dd\subset \ZZ$ to denote the set of negative fundamental discriminants, and we define:
\[\mathrm{supp}(\ell) = \set{(-d,m) \in \Dd \times \NN: (\exists a \in \ZZ)(a^2 - 4\ell = -dm^2)}\]
For every prime $\ell > 3$,
we have:
\[ \Delta_\ell(x) = x^{-2\left(1+\qr{-3}{\ell}   \right)}(x-1728)^{-\left(1+\qr{-4}{\ell}   \right)} \prod_{(-d,m)\in \mathrm{supp}(\ell)} H_{-d,m}(x)^{1 + \chi_\ell(-d)} \]
This factorization is discussed in detail in \cite{coxprimes},
see Theorem 13.2.

We will also use the following facts:
\begin{itemize}
    \item 
        Kronecker's congruence relation\footnote{See Ex. 2.20 in \cite{sil2}}
        implies that:
        \[ \Delta_\ell(x) \cong \prod_{j \in \FF_\ell} (x-j)^2 \pmod{\ell \ZZ[x]}\]
        This means that the factors $H_{-d,m}(x)$ of $\Delta_\ell$ all split over $\FF_\ell$.

    \item 
        By Theorem 14.16 in \cite{coxprimes}\footnote{See also Chapter 12 of \cite{langellfun}.}, the roots of $H_{-d,m}$ have the following interpretation in characteristic $\ell$:
        if $a^2 - 4\ell = (-d)m^2$, then the roots of $H_{-d,m}$ in $\FF_\ell$
        are precisely the $j$-invariants of elliptic curves which admit a model over $\FF_\ell$ with trace of Frobenius equal to $a$.
        In particular, the supersingular $j$-invariants in $\FF_p$
        are precisely the roots of $H_{-d,m}$, where $-dm^2 = -4\ell$.
    \item 
        
        For a discriminant $-d$, we write $h_{-d}(m)$ to denote the degree of $H_{-d,m}(x)$.
        Furthermore, we define a function $\varphi_{-d} : \NN\to \NN$:
        \[ \varphi_{-d}(m) = m \prod_{p|m} 1 - \frac{1}{p} \qr{-d}{p}\]
        Theorem 7.24 in \cite{coxprimes} shows that for $m > 1$, we have:
        \begin{equation}
            \label{eq:deghdm} h_{-d}(m) = \frac{2h_{-d}(1)}{|\Oo_d^\times|} \sum_{n | m} \phi_{-d}(n)
        \end{equation}
\end{itemize}

\subsection{Reduction mod $p$}

Let $p$ be a prime,
and let $\rho_{p,\ell} : X_0(\ell)_\Delta (\overline{\QQ})\to  X_0(\ell)_\Delta (\overline{\FF_p^{sep}})$ be the reduction map.
For a fixed value of $\ell$,
$\rho_{p,\ell}$ is surjective for all $p$ by Deuring lifting,
but $\rho_{p,\ell}$ can fail to be injective for certain primes $p$.
In order to use $\Delta_\ell$ in the computation of $\Gamma_{p,\ell}$,
we need to understand precisely when this happens.

{\lem{
Let $-d_1, -d_2$
and suppose $H_{-d_i,m_i}(x) | \Delta_\ell(x)$.

Let $P_i \in X_0(\ell)[H_{-d_i,m_i}]_\Delta$,
and let $\tilde{P_i}$ be the image of $P_i$ in $X_0(\ell)(\FF_p^{sep})$.

If $\tilde{P_1} = \tilde{P_2}$,
then $-d_1 = -d_2$.

}}

\begin{proof}
    Fix models $E_1,E_2$ in characteristic 0, and choose $E_1, E_2$
    with good reduction at some prime $\pr$ above $p$,
    and let $\xi_i \in End(E_i)$ be an endomorphism of degree $\ell$.
    Then $\xi_i$ satisfies a polynomial $x^2 - a_i x + \ell$,
    where $a_i^2 - 4\ell = (-d_i)m_i^2$.
    Note that $a_i$ is the trace of the associated endomorphism
    of the Tate module $T_\ell(E_i)$.
    
    Let $\tilde{E_i}$ be the reduction of $E_i$ mod $\pr$,
    and let $\tilde{\xi_i}$ be the reduction of $\xi_i$.
    Now, reducing mod $\pr$ induces an isomorphism
    of Tate modules $T_\ell(E_i) \to T_\ell(\tilde{E_i})$,
    so $\tilde{\xi_i}$ satisfies the same quadratic polynomial as $\xi_i$.

    Now, if $(\tilde{E_1},\ker(\tilde{\xi_1})),(\tilde{E_2},\ker(\tilde{\xi_2}))$
    are isomorphic, then the associated quadratics must split over the same extension of $\QQ$, so $-d_1 = -d_2$.

\end{proof}

Note that if $j(E_1), j(E_2)\not\in\set{0,1728}$,
then we can also conclude that $m_1 = m_2$ under the hypotheses of the previous lemma.
In any case, this allows to work with one factor $H_{-d,m}$ at a time.

Let $\Gg_d = Gal(L/\QQ)$ be as in \ref{sssec:notation},
and let $\FF_q$ be the splitting field of $(x^2+d)H_{-d}(x)$ over $\FF_p$.
Reduction mod $p$ induces a surjective homomorphism $\rho_\Gg : \Gg_d \to Gal(\FF_q/\FF_p)$,
and for $\tau \in \Gg_d$, $P \in X_0(\ell)[H_{-d,m}]$,
we have $\rho_{p,\ell}(\tau \cdot P) = \rho_\Gg(\tau) \cdot \rho_{p,\ell}(P)$.

{\prp{
Let $P \in X_0(\ell)[H_{-d,m}]$,
and assume $P \neq \mathrm{Fr}_\ell(P)$.

Then $\rho_{p,\ell}(P) = \mathrm{Fr}_\ell(\rho_{p,\ell}(P) )$
if and only if one of the following is true:
\begin{itemize}
    \item $\mathrm{Fr}_\ell(P) = P$.
    \item $\si_d \in \ker(\rho_\Gg)$.
\end{itemize}

}}.

\begin{proof}
    By Prop \ref{prp:frickechar0}, we have:
    \begin{align*}
        \mathrm{Fr}_\ell (\rho_{p,\ell}(P)) &= \rho_{p,\ell}(\mathrm{Fr}_\ell (P)) \\
        &= \rho_{p,\ell}(\si_d \cdot P) \\
        &= \rho(\si_d) \rho_{p,\ell}(P)
    \end{align*}
    If $\si_d \in \ker(\rho)$,
    then $\rho_{p,\ell}(P)$ is clearly a fixed point of $\mathrm{Fr}_\ell$ in characteristic $p$.
    Furthermore, if $P$ is fixed by Fricke in characteristic 0,
    then $\rho_{p,\ell}(P)$ is fixed by Fricke in characteristic $p$,
    which proves one direction.

    To prove the other direction,
    assume $\rho_{p,\ell}(P)$ is fixed by $\mathrm{Fr}_\ell$ and $\si_d \not\in \ker(\rho_\Gg)$.
    Then $\rho_\Gg(\si_d)$ is nontrivial, and it fixes $\rho_{p,\ell}(P)$,
    so $\si_d$ must fix $P$ in characteristic 0.
    By \ref{prp:frickechar0}, this is equivalent to $P$ being a fixed point of $\mathrm{Fr}_\ell$.
    
\end{proof}

This means that when we pass to characteristic $p$,
the number of fixed points of $X_0(\ell)$ under the Fricke involution
can go up, which means $\rho_{p,\ell}$ is not surjective.
However, this happens if and only if $\si_d \in \ker(\rho_\Gg)$.
Furthermore, note that when $H_{-d}$ has supersingular reduction, then $\qr{-d}{p} \neq 1$,
so $\si_d \in \ker(\rho_\Gg)$ is equivalent to $p | d$.

This explains what we observed in \ref{sssec:exp7l2}:
$\Delta_2(x)$ is divisible by $(x+3375)^2$ because the curve with CM by $\ZZ[\frac{1+\sqrt{-7}}{2}]$ gives rise to two points $P, \si(P)$ in characteristic 0,
but 7 ramifies in that ring, so $P, \si(P)$ map to the same point in $X_0(\ell)(\FF_7^{sep})$.

We will show that there is precisely one other reason that $\rho_{p,\ell}$
can fail to be injective: this happens when there exists $(-d,m) \in \mathrm{supp}(\ell)$ with $p|m$.
The smallest example of this phenomenon, if we exclude cases where $p\in \set{2,3}$,
is $p = 5, \ell =19$.

{\prp{
Let $p,\ell$ be primes,
and suppose there exists $(-d,pm) \in\mathrm{supp}(\ell)$
with $\qr{-d}{p} \neq 1$.

Then $X_0(\ell)[H_{-d,pm}]_\Delta(\FF_{p^2}) = X_0(\ell)[H_{-d,m}]_\Delta(\FF_{p^2})$.
}}

\begin{proof}
    Let $P' \in X_0(\ell[H_{-d,mp}(x)]_\Delta )(\overline{\QQ})$
    and assume $H_{-d,m}(j_\ell(P')) \neq 0$.
    Let $E'$ be an elliptic curve with $j(E') = j_\ell(P')$,
    and let $\xi' \in End(E')$ be an endomorphism of degree $\ell$
    such that $P' = [(E', \ker(\xi'))]$.

    By definition of $H_{-d,mp}$,
    we can find an isogeny $\psi : E \to E'$ of degree $p$
    from an elliptic curve $E$ that satisfies $H_{-d,m}(j(E)) = 0$,
    and $\xi \in End(E)$ such that 
    $\psi \circ \xi = \xi' \circ \psi$.
    Let $P = [(E, \ker(\xi))] \in X_0(\ell)(\overline{\QQ})$.

    Now, let $\tilde{E},\tilde{E'}$ be the reductions of $E,E'$ mod $p$.
    Note that $\tilde{E},\tilde{E'}$ are supersingular,
    so we may assume they are defined over $\FF_{p^2}$.

    The reduction $\tilde{\psi}$ is purely inseparable,
    so it factors as $\tilde{E}\to \tilde{E}^{(p)} \to \tilde{E'}$,
    where $\phi : \tilde{E}\to \tilde{E}^{(p)}$ is the map $\phi(x,y) = (x^p,y^p)$
    and where $\tilde{E}^{(p)} \to \tilde{E'}$ is an isogeny of degree 1.
    The second map is an isomorphism, 
    so we may assume that $\tilde{E'} = \tilde{E}^{(p)}$,
    and $\tilde{\psi} = \phi : (x,y)\mapsto (x^p,y^p)$.

    This means that $\tilde{P'} = \phi(\tilde{P}) \in X_0(\ell)[H_{-d,m}]_\Delta(\FF_{p^2})$.

\end{proof}

Finally, we show that these are the only reasons that $\rho_{p,\ell}$
can fail to be injective.\footnote{Note that this is true if we only care about points that map to $j_\ell\inv(S_p)$.}

{\lem{\label{lem:critsamej}
Let $L$ be a number field,
$E,E'/L$ elliptic curves and let $\psi : E \to E'$ be an isogeny,
with $\ker(\psi) \subset E(L)$.

Let $p$ be a prime and assume:
\begin{itemize}
    \item $E, E'$ have good reduction at a prime $\pr$ above $p$.
    \item $p\not| \deg(\psi)$.
\end{itemize}

Let $\tilde{\psi}: \tilde{E} \to \tilde{E'}$ be the reduction of $\psi$,
and let $End(\tilde{E})_0$ be the image of $End(E)$ in $End(\tilde{E})$.

If $\tilde{\psi} \in End(\tilde{E})_0$,
then $j(E) = j(E')$.

}}
\begin{proof}
    First, note that we may assume that $\ker(\psi)$ is cyclic:
    if $\ZZ/m_0 \times \ZZ/m_0 \subset \ker(\psi)$,
    then $\psi$ factors as $\psi = [m_0] \circ \psi_0$
    so we can replace $\psi$ by $\psi_0$ to obtain an isogeny with cyclic kernel.

    Now, let $H = \ker(\psi)$,
    and $\tilde{H}$ the reduction of $H$.
    Note that $\tilde{H}$ also has order $m$, since the reduction is injective on $m$-torsion for $m$ prime to $p$.
    Thus, if $P \in X_0(m)$ represents $(E,H)$,
    and $p\not| m$,
    then $\tilde{P} = [(\tilde{E},\tilde{H})]$ gives a point on $X_0(m)(\FF_p^{sep})$.

    If $\tilde{\psi} \in End(\tilde{E})_0$,
    that means $\tilde{E'} = \tilde{E}$ so $\tilde{P} \in X_0(m)_\Delta(\FF_p^{sep})$,
    and by Deuring lifting,
    we can lift $\tilde{P}$ to a point $P' \in X_0(m)_\Delta(\overline{\QQ})$,
    where $P'$ represents $(E,\ker(\xi))$ for some $\xi \in End(E)$.
    The kernels $\ker(\xi),\ker(\psi)$ are subgroups of $E$ of order $m$,
    and the reduction map is injective on $m$-torsion,
    so this implies that $\ker(\xi) = \ker(\psi)$ in characteristic 0.
    Thus, $E \cong E'$ and up to an isomorphism, we have $\xi = \psi$.

\end{proof}

{\prp{
Let $p, \ell$ be primes, let $-d < 0$ be a fundamental discriminant,
let $m > 0$ an integer and assume:
\begin{itemize}
    \item $\qr{-d}{p} \neq 1$
    \item $p \not| m$
\end{itemize}

Let $P,P' \in X_0(\ell)[H_{-d,m}]_\Delta$
and let $\tilde{P},\tilde{P'} \in X_0(\ell)(\FF_{p^2})$
be the images of $P,P'$ under the reduction map.

If $\tilde{P} = \tilde{P'}$, then $j_\ell(P) = j_\ell(P')$.

}}

\begin{proof}
    Fix representatives $(E,H), (E',H')$ of $P,P'$ in characteristic 0.
    Let $\psi_1 : E_0 \to E$ and $\psi_2 : E_0' \to E'$ be isogenies of degree dividing $m$ from elliptic curves $E_0, E_0'$ with CM by $\Oo_d$,
    and let $\psi_0 : E_0 \to E_0'$ be an isogeny of degree not divisible by $p$.

    Let $\psi : E\to E'$ be the composition $\psi = \psi_2 \circ \psi_0 \circ \widehat{\psi_1}$.

    Let $\xi \in End(E)$ be an endomorphism with $\ker(\xi) = H$.

    {\bf Claim 1:} $\xi \in \overline{End(E)_\psi}$.

    \begin{itemize}
        \item First, note that there exists $\xi_1 \in End(E_0)$ such that $\psi_1\circ \xi_1 = \xi \circ \psi_1$.
        \item Since $E_1,E_2$ have CM by $\Oo_d$,
        $F_{\psi_0} : End(E_1) \to End(E_2)$ is an isomorphism,
        so we can find $\xi_2 \in End(E_2)$ such that $(\psi_0 \circ \hat{\psi_1}) \circ \xi = \xi_2 \circ (\psi_0 \circ \hat{\psi_1})$.

        \item Finally, since $\deg \psi_2$ divides $m$,
        we can find $\xi' \in End(E')$ such that $\xi' \circ \psi = \psi \circ \xi$.
    \end{itemize}

    Let $H' = \ker(\xi') \subset E'$.

    {\bf Claim 2:} The pair $(E,H)$ is represented by $P'$ or $\mathrm{Fr}_\ell(P')$.

    This is because the ring $\Oo_d$ contains at most two ideals of norm $\ell$,
    and if there are two ideals, they are Galois conjugates.
    The Galois action $\sqrt{-d}\mapsto -\sqrt{-d}$ on $X_0(\ell)[H_{-d,m}]$
    coincides with the Fricke action by Prop \ref{prp:frickechar0}.
    Thus, replacing $P$ by $\mathrm{Fr}_\ell(P)$ if necessary,
    we may assume that $P'$ represents $(E,H')$.

    Suppose $\tilde{P} = \tilde{P'}$.
    Then $\tilde{E} \cong \tilde{E'}$,
    so composing with an isomorphism $\tilde{E'}\to \tilde{E}$ if necessary,
    we may assume that $\tilde{E'} = \tilde{E}$ and $\tilde{H'} = \tilde{H}$.

    This means $\tilde{\psi} \in End(\tilde{E})$,
    and $\tilde{\psi}$ commutes with $\tilde{\xi}$ in $End(\tilde{E})$.
    Now, let $End(\tilde{E})_0$ be the image of $End(E)\to End(\tilde{E})$.
    Note that $End(\tilde{E})_0\otimes \QQ$ is generated by $[1],\tilde{\xi}$,
    so in fact $\tilde{\psi}$ commutes with $End(\tilde{E})_0$.

    Finally, we use Lemma II.5.2 in \cite{sil2} to deduce that $\tilde{\psi} \in End(\tilde{E})_0$,
    and since $p\not|m$, 
    Lemma \ref{lem:critsamej} shows that $j(E) = j(E')$.
    
\end{proof}

\subsection{Main Results}

For a pair $-d,m$ and a prime $p$,
we define:
\[H_{-d,m,(p)} (x) = \frac{H_{-d}(x) H_{-d,m}(x)}{H_{-d,p^e}(x)}\]
where $p^e$ is the largest power of $p$ that divides $m$.
\begin{itemize}
    \item $H_{-d,m,(p)} \in \ZZ[x]$.
    \item $H_{-d,m,(p)}| H_{-d,m}$.
    \item $H_{-d,m,(p)} = H_{-d,m}$ if and only if $p\not| m$.
\end{itemize}

Define:
\[\Delta_{p,\ell}^{ss} = x^{-2\left(1+\qr{-3}{\ell}   \right)}(x-1728)^{-\left(1+\qr{-4}{\ell}   \right)} \prod_{(-d,m)\in \mathrm{supp}(\ell)} H_{-d,m,(p)}(x)^{\left(\frac{1-\qr{-d}{p}}{1+\chi_\ell(-d)}\right)} \]

This polynomial can be used to determine the number of endmorphisms of degree $\ell$ for a supersingular curve over $\FF_{p^2}$ in the same way we would use $\Delta_\ell$ in characteristic 0.

{\thm{
\label{thm:noedgesdeltap}
For $j_0 \in S_p$, the number of edges $j_0 \to j_0$  in $\Gamma_{p,\ell}$ 
is equal to the order of vanishing of vanishing of $\Delta_{p,\ell}^ss$ at $j_0$.

}}

Now, let $M_{p,\ell}$ be an adjacency matrix that represents $\Gamma_{p,\ell}$
relative to some ordering of the vertices, and write $Tr(\Gamma_{p,\ell})$ to denote the trace of $M_{p,\ell}$.
Then  $Tr(\Gamma_{p,\ell})$  is equal to the total number of edges $j_0 \to j_0$ in $\Gamma_{p,\ell}$.
Our analysis shows that $Tr(\Gamma_{p,\ell}) = \deg \Delta_{p,\ell}$.
This means we can obtain a formula for the trace that depends only on the degrees of $H_{-d,m,(p)}$.

Let $\ell>3$ be a prime, 
let $-d$ be a discriminant and define an integer $c_\ell(-d)$ as follows:
\begin{itemize}
    \item If $\ell$ does not split in $\Oo_d$, we set $c_\ell(-d) = 0$.
    \item If $-d \not\in \set{-3,-4}$ and $(-d,m) \in \mathrm{supp}(\ell)$,
    we define $c_\ell(-d) = \frac{2h_{-d}(m)}{1+\chi_\ell(-d)}$.
    \item If $\ell = a^2 + b^2$, where $0 < a < b$, we define $c_\ell(-4) = h_{-4}(ab)$.
    \item Finally, if $\ell = a^2 + 3b^2$,
    we define $c_\ell(-3) = h_{-3}(|b(a^2-b^2)|)$.
\end{itemize}

Then integer $c_\ell(-d)$ measures the contribution to $\Delta_\ell$
from elliptic curves whose endomorphism ring is isomorphic to an order in $\Oo_d$.
Since the degree of $\Delta_\ell$ is equal to $2\ell$,
and the contribution from any $-d$ is equal to $2c_\ell(-d)$,
it follows that $\sum_{-d \in \Dd} c_\ell(-d) = \ell$.

{\cor{
\label{cor:traceformula}

Let $p,\ell$ be primes, with $p \neq \ell$.

Then:
\[ Tr(\Gamma_{p,\ell}) = \ell + \left( \sum_{-d \in \Dd} \right) -\ep_{p,\ell}\]
where:
\[ \ep_{p,\ell} = \sum_{(-d,m) \in \mathrm{supp}(\ell)} \left(1 - \chi_\ell(d)\qr{-d}{p} \right)(h_{-d}(m) - h_{d}(m/p^e)) \]

Furthermore:
\begin{itemize}
    \item $\ep_{p,\ell} \geq 0$ for all $p,\ell$.
    \item For a fixed value of $\ell$, $\ep_{p,\ell} = 0$ for all but finitely many $p$.
\end{itemize}

}}

If we know $R_\ell(x,y)$, then we can compute $\Delta_\ell$
and read off the information needed to obtain the trace.
However, this information can also be obtained without needing any modular data:
the numbers $c_\ell(-d)$ can be interpreted as the number of elliptic curves $E/\FF_\ell$,
up to isomorphism over $\FF_\ell$,
for which the characteristic polynomial of Frobenius splits over $\Oo_d$.
Furthermore, using the formula \ref{eq:deghdm},
knowing $c_\ell(-d)$ is equivalent to knowing the class number of $\QQ(\sqrt{-d})$.
If we know $c_\ell(-d)$, or equivalently the degree of $H_{-d,m}$,
we can solve for the class number and then compute the degree of $H_{-d,m'}$
for every integer $m'$.
This means we can use the integers $c_\ell(-d)$ to obtain the value of $\ep_{p,\ell}$.

Altogether, this shows that for a fixed value of $\ell$,
we can obtain the value of $Tr(\Gamma_{p,\ell})$ for \emph{all} $p$ by counting points on elliptic curves over $\FF_\ell$,
see Algorithm \ref{alg:ecoverlcomp} for more details.

\subsection{Applications}

.

\subsubsection{$p \in \set{5,7,13}$}
\label{sssec:nakaya5713}
For $p \in \set{5,7,13}$,
there is precisely one supersingular $j$-invariant,
so every isogeny between supersingular curves is an endomorphism.
Computing the supersingular isogeny graph for these primes is trivial,
but, we can use our results to determine when $s_{p,\ell}^+$ will have
an irreducible quadratic factor, without needing to compute $s_{p,\ell}^+$.

The image of $P \in X_0(\ell)$ in $X_0(\ell)^+$ is defined over $\FF_p$
precisely when $\mathrm{Fr}_\ell (P) = \phi \cdot P$,
and $\phi \cdot P =\mathrm{Fr}_\ell(P)$ precisely when $\Hh_d \subset \ker(\rho_\Gg)$.

To give a concrete example, take $p = 13$ and $\ell = 5$.
We compute:
\[ \Delta_5(x) = H_{-20,1}(x) H_{-4,4}(x)^2H_{-11,1}(x)^2H_{-19,1}(x)^2\]

\begin{itemize}
    \item The factors $H_{-20,1}, H_{-11,1}, H_{-19,1}$ have supersingular reduction,
    but $-4 \equiv 1 \pmod 5$ is a nonzero quadratic residue,
    so we can ignore $H_{-4,4}$.
    \item $H_{-11}, H_{-19}$ have degree 1, so $\Hh_d$ is trivial. Thus, those endomorphisms give rise to points in $X_0(5)^+(\FF_{13})$.

    \item On the other hand, $H_{-20}$ has degree 2 and does not split over $\FF_{13}$,
    so $\Hh_5$ is not in the kernel of $\rho_\Gg$.
    Thus, the points that represent these isogenis should not be defined over $\FF_{13}$.
\end{itemize}

Altogether, this analysis suggests that $s_{13,5}^+$ has precisely one quadratic factor,
which is what we have in Table \ref{table:ob}.

We can also do this for $s_{5,13}^+$.
In this case, we have:
\[ \Delta_{p=5,\ell=13} = H_{-52,1}(x) H_{-3,12}(x)^2 H_{-43,1}(x)^2\]
Since $-3,-43$ have class number 1, they give rise to points in $X_0(13)^+(\FF_5)$,
but $H_{-52,1}(x)$ does not split over $\FF_5$, so the roots of that polynomial
give rise to a pair of conjugate points on $X_0(13)(\FF_{5^2})$.

\subsubsection{$p \in \set{11,17,19}$}
\label{sssec:tr11to19}

The set $S_p$ contains exactly 2 elements precisely when $p \in \set{11,17,19}$:
\[ S_{11} = \set{0,1} \qquad S_{17} = \set{0,3} \qquad S_{19} = \set{7,18} \]
Note that $1728 \equiv 18 \pmod{19}$
With the vertices ordered as above, the adjacency matrix of $\Gamma_{p,\ell}$ must have the form:
\begin{align*}
    M_{11,\ell} = \ttm{(\ell+1)-3a}{3a}{2a}{(\ell+1)-2a} & \quad p = 11 \\
    M_{17,\ell} = \ttm{(\ell+1)-3a}{3a}{a}{(\ell+1)-a} & \quad  p = 17 \\
    M_{19,\ell} = \ttm{(\ell+1)-a}{a}{2a}{(\ell+1)-2a} & \quad p = 19
\end{align*}
where $a$ is the number of points on $X_0(\ell)^+$ that represent isogenies between the two supersingular $j$-invariants.

If we know $Tr(\Gamma_{p,\ell})$, then we can solve for $a$
and obtain the complete graph.
In particular, this shows that $\Gamma_{p,\ell}$ is completely determined by the values of $h_{-d}(m)$ for $(-d,m) \in \mathrm{supp}(\ell)$ with $\qr{-d}{p} \neq 1$.

\subsubsection{$p \in \set{23,29,31}$}
If  $p \in \set{23,29,31}$,
the trace is no longer enough to determine the adjacency matrix,
since we have 3 parameters and the trace only allows us to eliminate one.
However, if we know the polynomials $H_{-d,m}(x)$,
then we can use Theorem \ref{thm:noedgesdeltap} to determine the number of endomorphisms $j_0\to j_0$
for each supersingular $j_0$.
This gives us three constraints, so we can solve for the adjacency matrix and thus obtain the graph.

Now, suppose $(-d,m) \in \mathrm{supp}(\ell)$ for all $\ell$ in some finite set $S$.
We can find the roots of $H_{-d,m}$ in $\FF_\ell$, and gives us the reduction of $H_{-d,m}$ in $\FF_\ell[x]$.
If $S$ is sufficiently large, then we can obtain a formula for $H_{-d,m}$ using the CRT method (\cite{Sutherland2009ComputingHC}).

\subsection{$p = 37$ and beyond}
\label{ssec:rosmc37}
Let $p = 37$.
Computing $\Gamma_{p,\ell}$ actually turns out to be easier in this case, for a few reasons.
\begin{itemize}
    \item First, note that $|S_p| = 3$, so we can compute $\Gamma_{p,\ell}$ using the polynomials $H_{-d,m}$ as we did with $23,29,31$.
    Thus, it is certainly no harder than the computation of $\Gamma_{p,\ell}$
    for those primes.

    \item Since, $p \equiv 1 \pmod{12}$, $j = 0,1728$ are not supersingular,
    so the adjacency matrix is actually symmetric in this case.

    \item Finally, $p = 37$ is \emph{not} supersingular.
    The supersingular $j$-invariants are $j_0 = 8$, and $j_1, j_2 \in \FF_{37^2}$,
    where $j_1, j_2$ are roots of $j^2 - 6j - 6 = 0$.

    We can use the Galois action to show that the number of edges $j_1 \to j_1$ is equal to the number of edges $j_2\to j_2$, and the number of edges $j_0 \to j_1$ is equal to the number of edges $j_0 \to j_2$.

\end{itemize}

If we could determine the number of edges $j_1 \to j_2$,
then we wouldn't need the actual polynomials $H_{-d,m}$ to compute $\Gamma_{p,\ell}$.
On a similar note,
if we could figure out which points $P\in X_0(\ell)$ satisfy $\si(P) = \mathrm{Fr}_\ell(P)$, then we could determine when points in $X_0(\ell)$
give rise to points on $X_0(\ell)^+(\FF_p)$.
This would allow us to extend the idea in \ref{sssec:nakaya5713} to $p \not\in \set{5,7,13}$.

To that end, we establish some notation.
Let $\ell$ be a prime, and let $Y_0(\ell)$ be the open subset of $X_0(\ell)$ obtained by removing the cusps of $X_0(\ell)$.
Note that $Y_0(\ell)$ is affine,
and the restriction of $j_\ell$ to $Y_0(\ell)$ gives a morphism $Y_0(\ell)\to \AA^1$.

Let $-d$ be a nonsquare in $\ZZ$ and let $K = \QQ(\sqrt{-d})$.
We write $Y_0(\ell)^{K/\QQ}$ to denote the restriction of scalars\footnote{See Ch 4.6 in \cite{poonenqpoints}} of $Y_0(\ell)\otimes \mathrm{Spec}\;K$:
\[ Y_0(\ell)_{K/\QQ} = \mathrm{Res}_{K/\QQ} (X_0(\ell) \times_{\QQ} \mathrm{Spec}\;K) \]

\begin{itemize}
    \item $Y_0(\ell)^{K/\QQ}$ is a 2-dimensional variety over $\QQ$.
    If $L/\QQ$ is a number field, then points on $Y_0(\ell)^{K/\QQ}(L)$ are in bijection with points on $Y_0(\ell)(K\otimes_\QQ L)$.

    \item The Galois group of $K/\QQ$ acts on $Y_0(\ell)_{K/\QQ}$,
    and the map $Y_0(\ell)^{K/\QQ} \to Y_0(\ell)^{K/\QQ}$ that takes $P$ to $\si_d(P)$ is an \emph{isomorphism of varieties}. 
    In contrast, the map $Y_0(\ell)(K) \to Y_0(\ell(K)$ defined in the same way is a bijection at the level of sets, but it is not a morphism of varieties in general.

\end{itemize}

Let $j_\ell^{K/\QQ} : Y_0(\ell)^{K/\QQ}  \to Y(1)^{K/\QQ}$ be the restriction of scalars of $j_\ell$,
and $\mathrm{Fr}_\ell^{K/\QQ}$ the restriction of scalars of the Fricke involution.
Note that $Y(1) \cong \AA^1$ so $Y(1)_{K/\QQ} \cong \AA^2$,
so we can add and subtract elements in $\AA^2$.
We use this to define a new function $\nabla : Y_0(\ell)_{K/\QQ}  \to Y(1)_{K/\QQ}$:
\[\nabla(P) = j_\ell^{K/\QQ} (\mathrm{Fr}_\ell^{K/\QQ} (P))- j_\ell^{K/\QQ}(\si_d(P)) \]

The points we seek lie on the subvariety:
\[ X_0(\ell)^{K/\QQ}_\nabla = \set{P \in X_0(\ell)^{K/\QQ}:
\nabla(P) = \mathbf{0}
}\]

Now, assume $X_0(\ell)$ has genus 0.
We can choose coordinates so that the cusps of $X_0(\ell)$ are at $0,\infty$.
Thus, $Y_0(\ell) \cong \mathrm{Spec} \; \QQ[t, t\inv]$.
Furthermore, the Fricke involution must have the form $\mathrm{Fr}_\ell(t) = \frac{n_\ell}{t}$ for some constant $n_\ell \in \QQ$.

Explicit models of this form for $\ell \in \set{3,5,7,13}$ can be found in the appendix of \cite{mestre}.
In all cases, the constant $n_\ell$ is an integer,
and in fact, $n_\ell$ is given by the simple formula $n_\ell = \ell^{e_\ell}$,
where $e_\ell = \lfloor \frac{18}{\ell} \rfloor$.

For an element $t \in K$ of norm $n_\ell$,
we have $\mathrm{Fr}_\ell(t) = \frac{t\si_d(t)}{t} = \si_d(t)$.
Thus, every element of norm $n_\ell$ in $K$ gives rise to a point on $X_0(\ell)^{K/\QQ}_\nabla$.
If $K$ contains any elements of norm $n_\ell$,
then it contains infinitely many such elements,
so $X_0(\ell)^{K/\QQ}_\nabla$ has positive dimension in this case.

There is a lot more one can say about $X_0(\ell)^{K/\QQ}_\nabla$,
but since our primary goal in this paper to describe $\Gamma_{p,\ell}$,
we will leave that to a future work. 

\section{Algorithms}

Let $p,\ell$ be primes, with $p\neq \ell$.
To compute the graph $\Gamma_{p,\ell}$,
we use the following meta-algorithm:

\begin{enumerate}
    \item First, we need to find \emph{any} supersingular $j$-invariant.
    We do this by looking for an elliptic curve $E_0/\FF_p$ with cardinality $p+1$.

    \item Next, we compute $\Gamma_{p,2}$ using models over $\FF_{p^2}$.
    This is computationally inexpensive and allows us to obtain the complete set $S_p$ starting from any single element.

    \item Finally, we use $X_0(\ell)^+$ to find all isogenies of degree $\ell$ between elements of $S_p$.
\end{enumerate}

\subsection{Step 1}

Our first step, regardless of which approach we choose to take,
is obtaining a model of \emph{any} supersingular curve.
This can be done without any sophisticated machinery, by taking advantage of the following two facts:

\begin{itemize}
    \item There is always at least one supersingular curve defined over $\FF_p$.
    \item We can determine whether an elliptic curve $E/\FF_p$ is supersingular by counting points on $E$: $E/\FF_p$ is supersingular if and only if $\# E(\FF_p) = p+1$.
\end{itemize}
Thus, we can get the process started by simply counting points on elliptic curves over $\FF_p$ until we find a curve with the desired cardinality.

There are various ways one can approach this search.
To simplify our next task, we will look for supersingular curves given by an equation of the form:
\begin{equation}
    y^2 = (x-c)(x^2 + ax + b) \qquad (a,b,c\in \FF_p)\label{eq:2tormodelgen}
\end{equation}

Setting $a = c = 0$ and $b = 1$,
we obtain a model for the curve with $j =1728$, which is supersingular whenever $p \equiv 3 \pmod 4$.
Thus, we may assume that $p \equiv 1 \pmod 4$.

This assumption implies that $\# E(\FF_p) = p+1 \equiv 2 \pmod 4$ for every supersingular curve $E/\FF_p$.
Thus every supersingular curve has precisely one nontrivial point of order 2,
so we deduce that:
\begin{itemize}
    \item Every supersingular curve over $\FF_p$ can be described by an equation of the form \ref{eq:2tormodelgen}.
    \item In any equation of that form which describes a supersingular curve, the quadratic factor is irreducible.
\end{itemize}

Fix a nonsquare $d \in \FF_p$,
and let $E/\FF_p$ be an elliptic curve given by an equation of the form:
\[ E: \quad y^2 = (x-c_0)(x^2+a_0x+b_0)\]
where the quadratic factor is irreducible.
We can complete the square to obtain a new equation of the form:
\[ E : \quad y^2 = (x-c_1)(x^2 - d_1)\]
where $d_1 = a^2 - 4b$ is the discriminant of the original quadratic.
Now, since the quadratic is irreducible, $d_1$ is a nonsquare in $\FF_p$,
so $d_1, d$ are in the same square class in $\FF_p^\times$,
i.e. there exists $t$ such that $t^2d_0 = d$.

\begin{itemize}
    \item If $t$ is a square - in other words, if $d,d_1$ are congruent modulo fourth powers - then we can use a change of variable of the form $(x,y)\mapsto (t^2 x, t^3y)$ to obtain an equation for $E$ of the form:
    \[  y^2 = (x-c)(x^2-d)\]
    \item If $t$ is not a square, then $E$ is not isomorphic to a curve given by an equation of that form over $\FF_p$, but the quadratic twist of $E$ admits a model of that form.
\end{itemize}

Altogether, this means that when $p \equiv 1 \pmod 4$,
we can find a representative for each supersingular $j$-invariant defined over $\FF_p$
in the family of curves:
\[ y^2 = (x-c)(x^2-d)\]
where $d$ is a fixed nonsquare in $\FF_p$,
and where $c$ ranges over all nonzero elements in $\FF_p$.

\begin{itemize}
    \item We ignore $c = 0$, as that always describes a curve with $j = 1728$, and 1728 is not supersingular when $p \equiv 1 \pmod 4$.
    \item We also only have to check one representative of each pair $\pm c$, since negating $c$ does not change the $j$-invariant.
\end{itemize}

Now, let $E$ be given by an equation as above.
The cardinality of $E$ is equal to:
\begin{align*}
    \# E(\FF_p) = 1+ \sum_{x \in \FF_p} \left( 1 + \left( \frac{(x-c)(x^2-d)}{p} \right) \right)
    = p+1 + \sum_{x\in \FF_p}    \left( \frac{x-c}{p}\right) \left( \frac{x^2-d}{p}\right)
\end{align*}
Thus $E$ is supersingular if and only if $\sum_{x\in \FF_p} \left( \frac{x-c}{p}\right) \left( \frac{x^2-d}{p}\right) = 0$.
This sum
is equal to the dot product of the vectors $v_c = \left( \qr{x-c}{p}\right)_{x \in \FF_p}, w_d = \left( \qr{x^2-d}{p}\right)){x\in \FF_p}$.
As $c$ varies, the vector $w_d$ does not change,
and the entries of $v_c$ simply gets rotated,
i.e. $(v_c)_x =(v_0)_{c+x}$.

\begin{algorithm}
\label{alg:sscount}
	\caption{Finding Supersingular Curve}
	\begin{algorithmic}
		\If{$p \equiv 3 \pmod 4$}
			\State \Return $(0,-1)$.
		\Else
			\State Compute $v_0 = \left( \qr{p}{x} \right)_{x \in \FF_p}$.
			\State Find nonsquare $d \in \FF_p^\times$.
			\State Compute $w_d = \left( \qr{p}{x^2-d} \right)_{x \in \FF_p}$
			\For{$\pm c \in \FF_p^\times$}
				\State Compute $v_c$ from $v_0$
				\State Compute $v_c \cdot w_d$.
				\If{$v_c \cdot w_d = 0$:}
					\State \Return $(2c,c^2-d)$.
			\EndIf
		\EndFor
		\EndIf
	\end{algorithmic}
\end{algorithm}

\subsection{Step 2}

In order to obtain the rest of the set $S_p$,
we will compute a supersingular isogeny graph using models.
The algorithm we use is well-known and is sometimes referred to as the \emph{graph method}, see e.g. \cite{mestre}, \cite{charleslauter}.

Let $\Ww$ be the space of Weierstrass equations over a field $\FF_p^{sep}$,
i.e. elements of $\Ww$ are 5-tuples $(a_1,a_2,a_3,a_4,a_6)$ that represent Weierstrass equations for elliptic curves.
For any $(a_i) \in \Ww$, and any $H \in S_\ell(E)$, where $E$ is the curve represented by $(a_i)$,
we can use V{\'e}lu's formulae (\cite{velu}) to obtain an explicit model in $\Ww$ that represents $E/H$.
Abstractly, this means there is a function $\Ww \to (\Ww^{\ell+1})/Sym(\ell)$
that takes $(a_i)$ to the unordered list of models that represent the quotients $E/H$.
The computation $\Gamma_{p,\ell}$ using models boils down to choosing an algorithm
$\nu_\ell$ that computes the image of $(a_i)$ in $\Ww^{\ell+1}$.

In general, we can do this using Algorithm \ref{alg:basistorsionsupsing}
but for $\ell = 2$, the whole process is much simpler.

Let $E$ be an elliptic curve given by an equation of the form:
\[ y^2 = x(x^2 + a_0x + b_0) \]
and let $r_1, r_2$ be the roots of $x^2+a_0 x + b_0 = 0$.
Then $E$ has 3 points of order 2, namely $P_0 = (0,0)$, $P_1 = (r_1, 0)$ and $P_2 = (r_2,0)$.
For $i = 1,2$, define $a_i = a_0 + 3r_i$ and $b_i = -a_0 r_i-2b_0$; these are the coefficients in the equation for $E$ where $P_1,P_2$ have been moved to the origin.
The isogenous curves are $E/\ip{P_i}$ are also described by equations of this form, with coefficients $(a_i',b_i') = (-2a_i,a_i^2-4b_i)$.
Thus, we can do all of our computations in the subset $\Ww_{a,b}$ instead of needing extending Weierstrass equations.
Furthermore, the hardest part in the computation of $\Ww_{a,b} \to \Ww_{a,b}^3/Sym(3)$ is solving the quadratic $x^2 + ax + b = 0$.

\subsubsection{Graph Method}

Let $p,\ell$ be primes, with $p\neq \ell$ and $p>3$,
let $n$ be the order of $p^2$ in $(\ZZ/\ell\ZZ)^\times$,
and let $q = p^n$.

To compute $\Gamma_{p,\ell}$ using the graph method,
we assume we are given a nonempty set \texttt{Models}
whose elements represent supersingular elliptic curves with full $\ell$-torsion over $\FF_q$.

\begin{algorithm}[ht!]
\label{alg:graphmethod}
	\caption{Graph Method}
	\begin{algorithmic}
            \Require{Set $\texttt{Models} \neq \emptyset$ as above.}
		\State Create an empty table $ \texttt{Graph} = \set{j : \set{}  \; j \in \emptyset}$.
		\While{$\mathrm{length}(\texttt{Models}) >0$} 
			\State Create empty list $\texttt{NewModels} = [ ]$.
			\For{$(a_i) \in \texttt{Models}$}
				\State Compute $j_0 = j((a_i))$.
				\If{$j_0 \not\in \texttt{Graph}$}
					\State{Compute $\nu_\ell(a_i)$.}
					\State{Update \texttt{Graph} by setting $\texttt{Graph}[j_0] = \set{j((a_i')) : (a_i') \in \nu_\ell((a_i))} $.}
					\State{Add the models in $\nu_\ell((a_i))$ to the list \texttt{NewModels}.}
				\EndIf
			\EndFor
			\State Set $\texttt{Models} = \texttt{NewModels}$.
		\EndWhile
		\State \Return \texttt{Graph}
	\end{algorithmic}
\end{algorithm}

The algorithm terminates when we stop finding new curves.
Since the $\ell$-isogeny graph is connected, we will stop finding new curves only
after we've seen every supersingular curve at least once,
so the output contains all of the information needed to construct $\Gamma_{p,\ell}$.
We can read off the set $S_p$ from the keys of \texttt{Graph}.

\subsection{Step 3}

Let $p,\ell$ be primes, and assume the following conditions are satisfied:
\begin{enumerate}
    \item We can compute the set of points $X_0(\ell)^+(\FF_{p^2})$.
    \item For each $Q = \set{P, \mathrm{Fr}_\ell(P)}\in X_0(\ell)^+(\FF_{p^2})$, we can determine the values of $a_\ell(Q) = j_\ell(P)+j_\ell(\mathrm{Fr}_\ell(P))$.
\end{enumerate}
For a supersingular prime $\ell$, both conditions are met if we a formula for the Atkin polynomial $R_\ell(x,y)$.

To obtain the graph, we use the following set:
\[ \Cc = \set{(j_1 + j_2, j_1j_2) : j_1, j_2 \in S_p, j_1 \neq j_2}\]
All we need to do is evaluate $a_\ell, b_\ell$ at every $Q \in X_0(\ell)^+$,
and record the roots of the quadratic if $(a_\ell(Q), b_\ell(Q)) \in \Cc$.

\begin{algorithm}
\label{alg:mcmethod}
	\caption{Modular Curve Method}
	\begin{algorithmic}
		\State Construct table $ \texttt{Counts} = \set{(a,b):0 |   \; (a,b)\in \Cc}$.
		\For{$Q \in X_0(\ell)^+(\FF_{p^2})$}
			\State Evaluate $a_\ell(Q), b_\ell(Q)$.
			\If{$(a_\ell(Q), b_\ell(Q)) \in \Cc$}
				\State Add 1 to $\texttt{Counts}[(a_\ell(Q),b_\ell(Q))]$.
			\EndIf
		\EndFor
		\State \Return \texttt{Counts}
	\end{algorithmic}
\end{algorithm}

By only looking for points that represent isogenies between distinct curves,
we don't have to worry about the issues involving endomrophisms.
Prop \ref{prp:x0plusprops} shows that there is no ambiguity regarding the number of edges represented by points which do not represent endomorphims.
Thus, we can use the output of this algorithm to construct $\Gamma_{p,\ell}^*$
and then obtain $\Gamma_{p,\ell}$ by adding the missing edges $j_0 \to j_0$.

\subsection{Trace Computations}

Let $\ell$ be a prime.
The following algorithm gives us the data needed to compute $Tr(\Gamma_{p,\ell})$
using Cor.\ref{cor:traceformula}.

\begin{algorithm}
\label{alg:ecoverlcomp}
	\caption{Trace Computation}
	\begin{algorithmic}
		\State Construct table $ \texttt{TraceFrobenius} = \set{a:\set{} |   \; 0 \leq a < 2\sqrt{\ell}}$.
            \State Compute the vectors $v_0 = (1)_{x \in \FF_\ell}, v_1 = (x)_{x \in \FF_\ell}$ and $v_3 = (x^3 \pmod \ell)_{\ell \in \FF_\ell}$.
		\For{$j_0 \in \FF_\ell \backslash \set{0,1728}$}
                \State Compute $t(j_0) = \frac{j_0}{1728-j_0} \in \FF_\ell$.
                \State Compute $v_{j_0} = v_3+3t(j_0)v_1 + 2t(j_0) v_0$.
                \State Compute $a_{j_0} = \left| \sum_{x \in \FF_\ell} \qr{(v_{j_0})_x}{-p} \right|$.
                \State Update $\texttt{TraceFrobenius}$ by adding $j_0$ to the $a_{j_0}$ entry.
		\EndFor
		\State \Return \texttt{TraceFrobenius}
	\end{algorithmic}
\end{algorithm}

The entry associated to each $a$ in the ouput
is the set of $j$-invariants that admit a model over $\FF_\ell$
with trace of Frobenius equal to $\pm a$.
We've left our the curves with $j = 0,1728$, but it is easy to add these back at the end:
\begin{itemize}
    \item If $p \equiv 2 \pmod 3$ (resp. $p \equiv 3 \pmod 4$), then $j = 0$ (resp. $j= 1728$) is supersingular, and therefore has trace 0.
    \item Otherwise, there are three (resp. two) values of $a > 0$ that give rise to a quadratic $x^2 - ax + \ell$ that splits over $\QQ(\sqrt{-3})$ (resp. $\QQ(i)$).
\end{itemize}

We can use the output of this algorithm to do the following:

\begin{enumerate}
    \item The $j$-invariants with trace of Frobenius $\pm a$ are the roots of $H_{-d,m}$ in $\FF_\ell$,
    where $a^2 - 4\ell = (-d)m^2$.
    Thus, we can determine the degree of $H_{-d,m}$, for all $(-d,m)\in \mathrm{supp}(\ell)$ directly from the output of this algorithm.
    \item We can solve for the class numbers of $\QQ(\sqrt{-d})$ using the degree of $H_{-d,m}$,
    and then compute the degree of $H_{-d,m'}$ for any integer $m'$.
    This means we can compute $Tr(\Gamma_{p,\ell})$ for any $p$.
    
    \item If we decide to compute a supersingular isogeny graph in \emph{characteristic} $\ell$ later on, then we already know the full set of supersingular $j$-invariants which are defined over $\FF_\ell$.

    \item If we run this algorithm for enough values of $\ell$, we can obtain formulae for (some) $H_{-d,m}$ using the CRT method \cite{Sutherland2009ComputingHC}.
    
\end{enumerate}

\subsection{Comparison}

For an integer $m$, we define $\ep(m) = 1$ if $m \equiv 2 \pmod 4$ and $\ep(m) = 0$
otherwise.
For $\ell$ is a prime and an integer $m$ that satisfies $\ell\not| m$, we write $\mathrm{ord}_{(\ZZ/\ell\ZZ)^\times} (m)$ to denote the order of $m$ in the multiplicative group $(\ZZ/\ell\ZZ)^\times$.
For a smooth curve $X$, write $\mathrm{gon}(X)$ to denote the gonality of $X$ and $\mathrm{gen}(X)$ to denote the genus of $X$.

We define:

\begin{align*}
    \delta_{GM}(p,\ell) &= 2\mathrm{ord}_{(\ZZ/\ell\ZZ)^\times} (p^2) - \ep(\mathrm{ord}_{(\ZZ/\ell\ZZ)^\times} (p)) \\
    \delta_{MC}(p,\ell)&= \mathrm{gon}(X_0(\ell)) \\
    \delta_{TR}(p,\ell)&= \mathrm{gen}(X_0(p)) 
\end{align*}

The integer $\delta_{GM}$, $\delta_{MC}$ and $\delta_{TR}$ quantify the difficulty
in computing $\Gamma_{p,\ell}$ using Algorithm \ref{alg:graphmethod}, Algorithm \ref{alg:mcmethod} and \ref{sssec:tr11to19}), respectively.

\begin{itemize}
    \item $\delta_{TR}(p,\ell) = 0$ if and only if $|S_p| = 1$.
    In this case, there is literally nothing to do to obtain the graph $\Gamma_{p,\ell}$ - it is simply the graph with 1 vertex and $\ell + 1$ edges.

    Note that $\delta_{GM}, \delta_{MC}$ are strictly positive for all $p,\ell$.

    \item The necessary computations can be done in $\FF_{p^2}$ or $\FF_\ell$ whenever $\delta(p,\ell) \leq 2$.

    \item In all 3 cases, the computations are a little simpler when $\delta(p,\ell) = 1$:
    \begin{itemize}
        \item In the graph method, the computation of $S_\ell(E_0)$ is easiest when $E_0$ is defined over $\FF_p$ and $\ell | p+1$.
        \item In the modular curve method, we have simple models of $X_0(\ell)$ in addition to the Atkin polynomials.
        \item In the trace method, we can obtain the complete graph from the integers $c_\ell(-d)$ without needing the polynomials $H_{-d,m}$.
        
    \end{itemize}
\end{itemize}

\section{Conclusion}
The algorithms we've described have been implemented in Python,
and the code can be found at the following repository:

\noindent \url{https://github.com/nhajouji/supsingecs}.

In addition, the repository contains precomputed tables that contain the following:
\begin{itemize}
    \item We ran Algorithm \ref{alg:ecoverlcomp} on every prime $\ell < 2^14$,
    and recorded the results in a table.
    This is the biggest table we computed, and it takes 58 MB to store as a pickled Pandas dataframe.
    
    \item We constructed a much smaller table that contain the degrees of $H_{-d,m}$ for all $(-d,m) \in \bigcup_{\ell < 2^{14}} \mathrm{supp}(\ell)$.

\end{itemize}

There are also additional tables, but we comment on how these tables can be used:
\begin{itemize}
    \item The second table requires much less space to store, and it contains all of the information needed to compute the trace.
    \item The first table is not redundant, however, since we can still use the more refined information about elliptic curves over $\FF_\ell$ in the CRT method.
    \item The $j$-invariants whose trace of Frobenius is equal to 0 are precisely the supersingular $j$-invariants in $\FF_\ell$.
    Thus, we can use part of the first table to start the process of computing $\Gamma_{p,\ell}$ using models, without needing Algorithm \ref{alg:sscount}.
\end{itemize}

\subsection{Generalizations and future work}
There are various ways one can try to extend the results in this paper.

\begin{enumerate}
    \item We've restricted our analysis of the algorithms to the case where $\delta(p,\ell)\leq 2$.
    The most obvious way of improving our work is to extend the analysis to settings where $\delta(p,\ell)>2$.

    \item The results in Section 3, as stated, only apply to supersingular curves.
    It should be possible to generalize the analysis to include ordinary curves, 
    although the statement may become more complicated.
    For example, we need a generalization of Prop \ref{prp:cmsupsingredprops} that includes the case where $\qr{-d}{p} = 1$.
    The results in \cite{hilbfac} should be helpful for this.

    \item A better understanding of the curves $X_0(\ell)^{K/\QQ}_\nabla$ would allow us to improve on the trace method for primes like $p = 37$,
    and it would allow us to further investigate the symmetry in Table \ref{table:ob} using arguments similar to \ref{sssec:nakaya5713}.

\end{enumerate}

\appendix

\section{Models of Modular Curves}
We briefly explain how to obtain the models needed for Algorithm \ref{alg:mcmethod}.

\subsection{Algorithm}

Let $\ell>3$ be a prime.
We can obtain models for $X_1(\ell),X_0(\ell),X_0(\ell)^+$ by doing the following:

\begin{enumerate}
	\item We construct the elliptic 3-fold $\Ee \to \AA^2 = \mathrm{Spec} \; \QQ[u,v]]$:
	\[ y^2 + (1-u) xy - v y = x^3 - vx^2 \]
	Let $j_\Ee : \AA^2 \to \PP^1$ be the $j$-invariant of $\Ee$.

	\item Let $P_0= (0,0) \in \Ee(\AA^2)$.
	We can compute $([m] P)= (x_m(u,v), y_m(u,v))$ for any integer $m \in \ZZ$.
	Note that $x_m, y_m \in \QQ(u,v)$.
	
	\item Let $a,b$ be integers and suppose $a+ b = \ell$.
	Setting $x([a]P_0) = x([b]P_0)$ gives us a relation $r_\ell(u,v)$ between $u,v$.
	Note that $r_\ell(u_0,v_0) = 0$ precisely when 
	the intersection of $P_0$ with the fiber over $(u_0,v_0)$ is a point of order $\ell$.

	\item Let $C_\ell \subset \AA^2$ be the curve where $r_\ell$ vanishes.
	For example, $r_5(u,v) = u-v$,
	and $r_7(u,v) = u^3 +uv-v^2$.
	Note that $C_5$ is smooth but $C_7$,
	is singular.
	
	Now, the curve $C_\ell$ is birational to $X_1(\ell)$,
	so we can obtain a model of $X_1(\ell)$ by computing the normalization map $X_1(\ell)\to C_\ell$.
	Furthermore, we define $X_1(\ell) \to \AA^2$ by composing the normalization map
	with the inclusion $C_\ell \to \AA^2$.
	This allows us to obtain a formula for the $j$-map $X_1(\ell) \to X(1)$ by composing the map $X_1(\ell)\to \AA^2$
	with $j_\Ee : \AA^2 \to X(1)$.

        The computations in \cite{x1nmodels} are very helpful for this step.
	
	\item Next, we obtain a model of $X_0(\ell)$ by computing the quotient of $X_1(\ell)$
	by the automorphism $\al: [(E,P)]\mapsto [(E,[m]P)]$,
	where $m$ is an integer whose image in $(\ZZ/\ell\ZZ)^\times/\ip{\pm 1}$ has order $\frac{\ell-1}{2}$.
        The results in \cite{autx13} can be used to obtain a formula for $\al$.
	
	Once we've computed the quotient map $X_1(\ell)\to X_0(\ell)$, we can obtain a formula for $j_\ell : X_0(\ell) \to X(1)$ by finding a function on $X_0(\ell)$
	that pulls back to the correct function $X_1(\ell)\to X(1)$.
	
	\item We also need a formula for the Fricke involution.
	This can be obtained in several ways:
	\begin{itemize}
		\item We can apply V{\'e}lu's formula to the whole family $\Ee|_{X_1(\ell)} \to X_1(\ell)$,
		and obtain a formula for $j_\ell \circ \mathrm{Fr}_\ell$ the same way we obtained $j_\ell$.
		This always works, but the equations get increasingly complicated as $\ell$ increases.
		
		\item 
		When the genus of $X_0(\ell)$ is bigger than 1,
		the automorphism group of $X_0(\ell)$ is finite, and the Fricke involution is one of the elements in that finite group.
		The Fricke involution swaps the cusps of $X_0(\ell)$, so we can rule out any involutions that don't swap the cusps.
            In certain cases, this may be enough to determine $\mathrm{Fr}_\ell$.
	\end{itemize}
	
	\item Finally, once we have a formula for $j_\ell$ and $\mathrm{Fr}_\ell$,
	computing the quotient $X_0(\ell)\to X_0(\ell)^+$ and the functions $a_\ell, b_\ell$ is straightforward.

\end{enumerate}

To give a concrete example,
we use this to obtain models for the modular curves when $\ell = 11$.
First, note that $X_1(\ell)$ is isomorphic to the elliptic curve:
\[ t^2 - t = s^3 - s^2 \]
and the map $X_1(\ell) \to \AA^2$ is given by:
\[(s,t) \mapsto (u,v) = \left( \frac{(t-1) (s+t-1)}{s},\frac{t(t-1)  (s+t-1)}{s}\right) \]

\begin{itemize}
	\item The Mordell-Weil group of $X_1(11)$ over $\QQ$ has order 5,
	and in fact, this curve is isomorphic to the fiber over $(u,v) = (1,1) \in C_5 \subset \AA^2$.
	The Mordell-Weil group can be generated using the point $(0,0)$.
	
	\item The automorphism group of $X_1(11)$ is generated by the translation maps and the negation map $(s,t)\mapsto (s, s-t)$.
	We deduce that $Aut(X_1(11)/X_0(11))$ is generated by the translation map $P\mapsto (0,0)+P$,
	which means that the map $X_1(11) \to X_0(11)$ is the isogeny of degree 5 whose kernel contains $(0,0)$.
	We can use V{\'e}lu's formula to obtain an explicit formula for the isogeny $X_1(\ell)\to X_0(\ell)$,
	as well as a model for the codomain:
	\[w^2 -w = z^3 - z^2 - 10z - 20\]
	
	 \item One of the cusps of $X_0(11)$ lies at the point at infinity,
	 and the other cusp is at the point $(z_0,w_0) = (16,61)$.
	 The only involution that swaps these points is $P\mapsto (16,61)-P$,
	 where the difference is computed using the group law on $X_0(11)$,
	 so we deduce that $\mathrm{Fr}_{11}(P) = P\mapsto (16,61)-P$.
	 We can also obtain a formula for the $j$-map:
	\[j_{11}((z,w)) = \frac{\left(z^4+12 z^3-8 z^2 w-114 z^2+40 z w+140 z-32 w+217\right)^3}{(-5 z+w+19) \left(4 z^2-z w-29 z+4 w+51\right)^2} \]
	 
	 \item Let $z : X_0(11)\to \PP^1$ be the function $z((z_0,w_0)) = z_0$
	 and let
	  $x : X_0(11)\to \PP^1$ be the function $x(P) = z(P+[3](16,61))$.
	  
	  Then $x$ is a function of degree 2, and $x$ is invariant under the Fricke involution,
	  so $x$ generates the subfield of $K(X_0(11))$ fixed by $\mathrm{Fr}_{11}$.
	  Note that this subfield is isomorphic to the function field of $X_0(11)^+$.
	  Furthermore, $P$ is a cusp if and only if $x(P) = 5$.
	  
	  \item The functions $j_{11}+j_{11} \circ \mathrm{Fr}_{11},j_{11}\cdot j_{11} \circ \mathrm{Fr}_{11}$ are
	   rational functions in $x$,
	   and after applying a suitable fractional linear transformation,
	   we can simply read off the coefficients of the polynomials $a_{11},b_{11}$.

\end{itemize}

We can do this for any $\ell$, but the computations and functions involved get increasingly complex
as the genus of $X_1(\ell)$ goes up.

\subsection{Models of $X_0(\ell)$}
\label{ssec:modelsx0}

For $\ell \in \set{3,5,7,11,13}$, the curve $X_0(\ell)$ has genus 0.
Formulas for the maps $j_\ell, \mathrm{Fr}_\ell$ can be found at the end of \cite{mestre}.
We record them here, since they are needed to verify the results in \ref{ssec:rosmc37}.

\begin{align*}
j_3&= \frac{(27-t)\left.(t-3)^3\right.}{t} & & \mathrm{Fr}_3(t) =\frac{729}{t} = \frac{3^6}{t}\\
j_5 &=\frac{\left(t^2+10 t+5\right)^3}{t} & & \mathrm{Fr}_5(t) =\frac{125}{t} = \frac{5^3}{t}\\
j_7 &=\frac{-\left(t^2-13 t+49\right) \left(t^2-5 t+1\right)^3}{t} & & \mathrm{Fr}_7(t) =\frac{49}{t} = \frac{7^2}{t}\\
j_{13} &=\frac{-\left(t^2-5 t+13\right) \left(t^4-7 t^3+20 t^2-19 t+1\right)^3}{t} & & \mathrm{Fr}_{13}(t) =\frac{13}{t}\\
\end{align*}

\subsection{Models of $X_0(\ell)^+$}
\label{ssec:modelsatkin}

Formulae for the Atkin polynomials $R_2, R_3, R_5, R_7$ can be found in several places,
e.g. \cite{x0plusconjecture},
but the formulae for $\ell > 7$ are harder to find.
The formula for $\ell = 11$ that we just obtained is:

\begin{align*}
    R_{11}(x,y) &= x^2 - ( 
     y^{11}-44 y^{10}+693 y^9-4334 y^8+4400 y^7 +42658 y^6\cdots \\
    & -44968 y^5-178376 y^4-58432 y^3+86240 y^2+67200 y+16000)x +\\
     &\left(y^4+232 y^3+1176 y^2+1120 y+400\right)^3
\end{align*}

Fortunately, SageMath\cite{sagemath} contains formulae for $R_\ell$ for all supersingular primes $\ell$.
\footnote{Precisely, the formula can be found in one of several databases available in SageMath. There is also a database containing formulae for the classical modular polynomials.}
For example, we found the following formulae for $R_{17}, R_{19}$:

\begin{equation*}
\begin{split}
R_{17}(x,y)&=x^2 - (y^{17}-119 y^{15}-238 y^{14}+5338 y^{13}+20808 y^{12}-91766 y^{11}-\cdots \\
&630836 y^{10}-70737 y^9+7118240 y^8+16023299 y^7-13049914 y^6-98725154 y^5-\cdots \\
&125706976 y^4+18738794 y^3+169635044 y^2+128884056 y+25608112)x+\\
&\left(y^6+248 y^5+4082 y^4+25988 y^3+80561 y^2+122444 y+73252\right)^3
\end{split}
\end{equation*}

\begin{equation*}
\begin{split}
R_{19}(x,y)&=x^2 -(y+3) \left(y^6-y^5-36 y^4-30 y^3+246 y^2+432 y+160\right) \cdot \\
&(y^{12}-2 y^{11}-71 y^{10}+12 y^9+1848 y^8+2532 y^7-17359 y^6-46922 y^5+\cdots \\
&12297 y^4+145008 y^3+119712 y^2-33600 y-51200)x + \cdot\\
&(y+3)^2 \left(y^6+246 y^5+3593 y^4+19312 y^3+48544 y^2+57280 y+25600\right)^3
\end{split}
\end{equation*}

\bibliographystyle{plain}
\bibliography{refs}

\end{document}